\pdfoutput=1
\documentclass[reqno]{amsart}
\usepackage{enumerate}
\usepackage{tabto}
\usepackage[mathscr]{euscript}
\usepackage{xcolor}
\usepackage{layout}
\usepackage{fancyhdr}
\usepackage{array}
\usepackage{amsfonts}
\usepackage{amsmath}
\usepackage{amssymb}
\usepackage{mathtools}
\usepackage{graphicx}
\usepackage{bm}
\usepackage{enumitem}
\usepackage{caption}
\usepackage{color}
\usepackage{csquotes}
\usepackage{bookmark}
\usepackage{float}
\usepackage{multirow}
\usepackage[square,numbers,sort&compress]{natbib}
\usepackage{hyperref}
\hypersetup{colorlinks=true,linkcolor=blue,citecolor=red}
\allowdisplaybreaks

\def\XXint#1#2#3{{\setbox0=\hbox{$#1{#2#3}{\int}$ }
\vcenter{\hbox{$#2#3$ }}\kern-.6\wd0}}

\newtheorem{theorem}{Theorem}[section]
\newtheorem{lemma}[theorem]{Lemma}
\newtheorem{corollary}[theorem]{Corollary}
\newtheorem{proposition}[theorem]{Proposition}
\newtheorem{remark}[theorem]{Remark}
\theoremstyle{definition}
\newtheorem{definition}[theorem]{Definition}

\numberwithin{equation}{section}

\newcommand{ \mr }{ \mathbb{R} }

\newcommand{\iintss}{\mathop{\mathchoice
  {\displaystyle\int\mkern-14mu\int}
  {\textstyle\int\mkern-7mu\int}
  {\scriptstyle\int\mkern-7mu\int}
  {\scriptscriptstyle\int\mkern-7mu\int}}\nolimits}
\newcommand{\iints}[1]{\iintss_{\!#1}}
\newcommand{\miints}{\mathop{\vphantom{\iint}\mathchoice
  {\ooalign{\hidewidth$\displaystyle-\mkern-6mu-$\hidewidth\cr$\displaystyle\int\mkern-14mu\int$\cr}}
  {\ooalign{\hidewidth$\textstyle-\mkern-6mu-$\hidewidth\cr$\textstyle\int\mkern-7mu\int$\cr}}
  {\ooalign{\hidewidth$\scriptstyle-\mkern-6mu-$\hidewidth\cr$\scriptstyle\int\mkern-7mu\int$\cr}}
  {\ooalign{\hidewidth$\scriptscriptstyle-\mkern-6mu-$\hidewidth\cr$\scriptscriptstyle\int\mkern-7mu\int$\cr}}}\nolimits}
\newcommand{\miint}[1]{\miints_{\!#1}}

\newcommand{\osc}{\operatorname*{osc}}
\newcommand{\loc}{\operatorname{loc}}

\begin{document}
\title[Gradient regularity without gap conditions]{Gradient regularity without gap conditions for parabolic double phase systems with coefficients depending only on time}

\author{Jehan Oh}\address{Department of Mathematics, Kyungpook National University, Daegu, 41566, Republic of Korea} \email{jehan.oh@knu.ac.kr}


\begin{abstract}
We study degenerate parabolic double phase systems with growth exponents $p<q$ whose modulating coefficient depends only on time. In this case, the operator on each time slice has Orlicz growth with structural constants independent of the coefficient, and the system is invariant under spatial translations. Using these features, we prove second order estimates and local bounds for the gradient when the coefficient is merely bounded and measurable, and $C^{1,\alpha}$ regularity when it is H\"older continuous in time. For weak solutions with gradients in $L^q$, these results hold without any gap condition between $p$ and $q$.
\end{abstract}
\maketitle

\section{\bf Introduction}\label{section 1}
In this paper we study the regularity of the spatial gradient of weak solutions to parabolic double phase systems of the form
\begin{equation}    \label{eq : main equation}
    u_t - \operatorname{div} \mathcal{A}(t,Du)=0 \qquad \text{in} \ \, \Omega_T\coloneq \Omega\times (0,T),
\end{equation}
where
\begin{equation}    \label{eq : vector field}
    \mathcal{A}(t,\xi)\coloneq |\xi|^{p-2}\xi+a(t)|\xi|^{q-2}\xi \qquad \text{for } t\in(0,T) \text{ and } \xi\in\mr^{N\times n}.
\end{equation}
Here $n\geq2$, $\Omega\subset\mr^n$ is a bounded open set, $u:\Omega_T\to\mr^N$ with $N\geq1$, the growth exponents satisfy $2\leq p<q<\infty$, and the modulating coefficient $a:(0,T)\to[0,\infty)$ is a bounded measurable function which depends only on the time variable. We write
$$
H(t,\varkappa)\coloneq \varkappa^p+a(t)\varkappa^q \qquad \text{for } \varkappa\geq0 \text{ and } t\in(0,T),
$$
so that $\mathcal{A}(t,\xi)\cdot\xi=H(t,|\xi|)$. Our main results state that weak solutions to \eqref{eq : main equation} have locally bounded gradients when $a$ is merely bounded and measurable, and locally H\"older continuous gradients when $a$ is H\"older continuous in time, without any condition relating $p$ and $q$. In this sense the regularity theory of \eqref{eq : main equation} is gap-free, in contrast with the theory of double phase problems with coefficients depending on the space variable, which we now recall.

Parabolic double phase systems are the evolutionary counterpart of the double phase functional
$$
W^{1,1}(\Omega) \ni w \mapsto \int_\Omega \left[\frac{1}{p}|Dw|^p+\frac{1}{q}a(x)|Dw|^q\right] dx,
$$
which was introduced by Zhikov \cite{Zhikov1986,Zhikov1993,Zhikov1995,Zhikov1997} in the study of homogenization and of the Lavrentiev phenomenon, and which belongs to the class of functionals with $(p,q)$-growth whose study was initiated by Marcellini \cite{Marcellini1989,Marcellini1991}. The functional models composite media whose hardening exponent changes from $p$ to $q$ across the zero set of the coefficient, see \cite{DeFilippisMingione2023}. For the elliptic problem, the regularity of minimizers is governed by the interplay between the difference $q-p$ and the regularity of the coefficient. If $a\in C^{0,\alpha}(\Omega)$ with $\alpha\in(0,1]$, local minimizers, bounded minimizers and H\"older continuous minimizers with exponent $\gamma_0$ have H\"older continuous gradients under the respective conditions
$$
\frac qp\leq1+\frac\alpha n,\qquad q\leq p+\alpha,\qquad q<p+\frac{\alpha}{1-\gamma_0},
$$
see \cite{Colombo2015,Baroni2018}, \cite{Colombo2015a,Baroni2018} and \cite{Baroni2018}, respectively. The first two conditions are sharp in view of the counterexamples in \cite{Esposito2004,Fonseca2004,Balci2020}, see also \cite{Mingione2021} for a survey.

In the parabolic setting the coefficient is usually allowed to depend on $z=(x,t)$ and to belong to the parabolic H\"older class $C^{\alpha,\frac{\alpha}{2}}(\Omega_T)$. Kim, Kinnunen and Moring \cite{2023_Gradient_Higher_Integrability_for_Degenerate_Parabolic_Double-Phase_Systems} proved gradient higher integrability for degenerate parabolic double phase systems under the gap bound
$$
q\leq p+\frac{2\alpha}{n+2},
$$
and the parabolic Lipschitz truncation of Kim, Kinnunen and S\"arki\"o \cite{Wontae2023a} provides energy estimates in the same range. These results were extended to the singular range and to Calder\'on--Zygmund type estimates \cite{Wontae2024,Wontae2023b} and to generalized Orlicz growth \cite{HastoOk2025}. For bounded solutions and for H\"older continuous solutions with exponent $\gamma_0$, the gap bound was relaxed to
$$
q\leq p+\alpha\qquad\text{and}\qquad q\leq p+\frac{\alpha}{1-\gamma_0}\quad\text{with}\quad q<p+1,
$$
respectively, see \cite{KimOh2026,Chlebicka2025} and \cite{Oh2026Holder}. We refer to \cite{Kim2024,Kim2025,Sen2025,SenSiltakoski2026} for multi-phase problems and improved gap bounds, to \cite{Buryachenko2022,Wontae2025} for the continuity and the H\"older continuity of bounded solutions, and to \cite{KimKimOh2026,KimOh2026b} for the absence of the Lavrentiev phenomenon. Concerning the regularity of the gradient, Sen and Siltakoski \cite{SenSiltakoski2025} used viscosity methods to prove Lipschitz continuity in space of bounded solutions when $q\leq p+1$ and the coefficient is Lipschitz continuous in space, continuous in time and almost nondecreasing in time. To the best of our knowledge, all these results require a gap bound, and, unlike in the elliptic case, it is not known whether any of the parabolic gap bounds is sharp.

Time-only coefficients, which are the subject of this paper, have no elliptic counterpart, since in a stationary problem they reduce to constants, and the double phase functional with a constant coefficient has autonomous Orlicz growth, for which no gap condition arises \cite{Lieberman1991}. In the parabolic setting they describe a medium whose second phase is present only at the times at which the coefficient is positive. Time-only coefficients are, however, included in the theory of B\"ogelein, Duzaar and Marcellini \cite{BDM2013} for parabolic equations with $(p,q)$-growth, where the vector field is differentiable in the spatial variables and only measurable in time. For weak solutions $u\in L^q_{\loc}(0,T;W^{1,q}_{\loc}(\Omega))$ to \eqref{eq : main equation} with $N=1$, their results give a local bound for the gradient and the existence of solutions to the Cauchy--Dirichlet problem under the conditions
$$
2\leq p\leq q<p+\frac4n\qquad\text{and}\qquad 2\leq p\leq q<p+\frac4{n+2},
$$
respectively, see \cite[Section 1.3]{BDM2013}. In \cite{BDM2013}, the bound on $q$ appears in a Moser iteration for the differentiated equation, whose ellipticity is of order $p$ while its growth is of order $q$. The subquadratic case was treated by Singer \cite{Singer2015}, and De Filippis \cite{DeFilippis2020} constructed solutions with locally bounded gradients also for coefficients which are Sobolev functions in space, under a stronger restriction on $q-p$. Without any growth condition, Sch\"atzler and Siltakoski \cite{SchatzlerSiltakoski2023} obtained, for time-dependent integrands, variational solutions which are Lipschitz continuous in space, provided that the domain is convex and the data satisfy the bounded slope condition.

Our starting point consists of two elementary observations. First, for every $t$ the vector field $\xi\mapsto\mathcal{A}(t,\xi)$ has the Uhlenbeck structure
$$
\mathcal{A}(t,\xi)=\mathcal{G}(t,|\xi|)\,\xi,\qquad \mathcal{G}(t,s)\coloneq s^{p-2}+a(t)s^{q-2},
$$
and
\begin{equation}\label{cond : Lieberman condition}
0\leq s\,\partial_s\mathcal{G}(t,s)\leq (q-2)\,\mathcal{G}(t,s)\qquad\text{for all } s>0,
\end{equation}
where the bounds do not depend on $t$ or on the value $a(t)\geq0$. Consequently the eigenvalues of $\partial_\xi\mathcal{A}(t,\xi)$ lie between $\mathcal{G}(t,|\xi|)$ and $(q-1)\mathcal{G}(t,|\xi|)$, and on every time slice the operator in \eqref{eq : main equation} belongs to the same class of operators with Orlicz growth in the sense of Lieberman \cite{Lieberman1991}, whether or not $a(t)$ vanishes. In the terminology of De Filippis and Mingione \cite{DeFilippisMingione2021,DeFilippisMingione2023}, the double phase operator is uniformly elliptic in the pointwise sense, since the ratio between the largest and the smallest eigenvalues of $\partial_\xi\mathcal{A}$ is bounded by $q-1$. Its nonlocal ellipticity ratio, in which the largest eigenvalue is maximized and the smallest one is minimized over a cylinder $Q$, is instead of order $1+\|a\|_{L^\infty(Q)}|\xi|^{q-p}$ whenever $\inf_Qa=0<\|a\|_{L^\infty(Q)}$, and hence unbounded in $\xi$. For time-only coefficients the nonlocal ellipticity ratio is still unbounded on space-time cylinders, but on every time slice it is bounded by $q-1$. Second, the system \eqref{eq : main equation} is invariant under translations in the spatial variables. When $a$ depends also on $x$, the first property holds only pointwise, since the growth function changes from point to point in each time slice, and the second one fails.

The gap bounds in the parabolic literature quantify the extent to which the pointwise uniform ellipticity can compensate for the unboundedness of the nonlocal ellipticity ratio on cylinders. In the intrinsic approach of \cite{2023_Gradient_Higher_Integrability_for_Degenerate_Parabolic_Double-Phase_Systems}, for instance, the gap bound comes from the requirement that the oscillation of the coefficient on intrinsic cylinders be absorbed by the $p$-phase. For time-only coefficients this oscillation is smaller, since intrinsic cylinders are short in time in the degenerate range, but it does not vanish, and a restriction on $q-p$ appears to remain. We show that for time-only coefficients no gap condition is needed at all. For weak solutions in the sense of Definition \ref{def : weak solution} below, we prove second order estimates, higher integrability with an explicit exponent, local bounds and H\"older continuity of the gradient for every $N\geq1$ and $2\leq p<q<\infty$. These appear to be the first regularity results for parabolic double phase systems with a nonconstant coefficient which hold for arbitrary weak solutions and involve no condition relating $p$ and $q$.

We now define weak solutions. Here and below, $z=(x,t)$ denotes a point of $\Omega_T$ and $dz=dx\,dt$.
\begin{definition}\label{def : weak solution}
    A function $u:\Omega_T\to\mr^N$ with
    $$
    u\in C(0,T;L^2(\Omega,\mr^N))\cap L^q(0,T;W^{1,q}(\Omega,\mr^N))
    $$
    is a weak solution to \eqref{eq : main equation} if
    $$
    \iints{\Omega_T} \left[-u\cdot \varphi_t +\mathcal{A}(t,Du)\cdot D\varphi\right] dz=0
    $$
    for every $\varphi\in C_0^\infty(\Omega_T,\mr^N)$.
\end{definition}

As in \cite{Oh2026Holder}, the preceding definition includes the $L^q$-integrability of the gradient, and a local version of this assumption is also used in \cite{BDM2013}. In this paper, the assumption is used only qualitatively, namely to justify testing the system with functions built from $u$ and to carry out the approximation argument of Section \ref{section 4}, and none of our estimates depends on the $L^q$-norm of $Du$. The natural energy class would instead require
$$
u\in C(0,T;L^2(\Omega,\mr^N))\cap L^p(0,T;W^{1,p}(\Omega,\mr^N)),\qquad
\iints{\Omega_T} H(t,|Du|)\,dz<\infty.
$$
Extending our results to this class requires a separate justification of the energy estimates, which is related to the absence of a Lavrentiev phenomenon in the time variable, and we do not pursue this here. We also note that the gap bounds in \cite{BDM2013,2023_Gradient_Higher_Integrability_for_Degenerate_Parabolic_Double-Phase_Systems,KimOh2026,Oh2026Holder} are needed in the proofs there even for solutions in the class of Definition \ref{def : weak solution}.

In what follows, $Q_r(z_0)\coloneq B_r(x_0)\times\Lambda_r(t_0)$ with $\Lambda_r(t_0)\coloneq(t_0-r^2,t_0+r^2)$ denotes the parabolic cylinder centered at $z_0=(x_0,t_0)$, and
$$
V_m(\xi)\coloneq|\xi|^{\frac{m-2}{2}}\xi\qquad\text{for } m\geq2 \text{ and } \xi\in\mr^{N\times n}.
$$
Our first result is a second order estimate for merely bounded and measurable coefficients. In its statement, $a(t)|DV_q(Du)|^2$ stands for $|D(a^{\frac12}V_q(Du))|^2$, where $a^{\frac12}V_q(Du)$ denotes the function $(x,t)\mapsto a(t)^{\frac12}V_q(Du(x,t))$.
\begin{theorem}\label{thm : second order}
    Let $a\in L^\infty(0,T)$ be nonnegative and let $u$ be a weak solution to \eqref{eq : main equation}. Then, for every $Q_{2r}(z_0)\Subset\Omega_T$, the functions $V_p(Du)$ and $a^{\frac12}V_q(Du)$ belong to
    $$
    L^2(\Lambda_r(t_0);W^{1,2}(B_r(x_0),\mr^{N\times n})),
    $$
    and
    \begin{equation}\label{eq : second order estimate}
    \begin{aligned}
        &\sup_{t\in\Lambda_r(t_0)}\int_{B_r(x_0)}|Du|^2\,dx+\iints{Q_r(z_0)}\left[|DV_p(Du)|^2+a(t)|DV_q(Du)|^2\right] dz\\
        &\qquad\leq \frac{c}{r^2}\iints{Q_{2r}(z_0)}\left[|Du|^2+H(t,|Du|)\right] dz
    \end{aligned}
    \end{equation}
    for some constant $c=c(n,p,q)$.
\end{theorem}

Theorem \ref{thm : second order} extends the second order estimate for the degenerate $p$-Laplace system, see \cite{DiBenedettoFriedman1984} and \cite[Chapter VIII]{1993_Degenerate_parabolic_equations_DiBenedetto}, to double phase systems with time-only coefficients, and its constant does not depend on the coefficient. Combined with the Sobolev inequality on time slices, it yields the higher integrability $H(t,|Du|)^{1+\frac{4}{nq}}\in L^1_{\loc}(\Omega_T)$ with an explicit exponent, see Corollary \ref{cor : higher integrability}.

The second result is a local bound for the gradient, that is, a local Lipschitz bound for $u$ in the spatial variables, again for bounded and measurable coefficients.
\begin{theorem}\label{thm : Lipschitz}
    Let $a\in L^\infty(0,T)$ be nonnegative and let $u$ be a weak solution to \eqref{eq : main equation}. Then $Du\in L^\infty_{\loc}(\Omega_T,\mr^{N\times n})$, and there exists a constant $c=c(n,q)\geq1$ such that
    \begin{equation}\label{eq : Lipschitz estimate}
        \sup_{Q_{r/2}(z_0)}|Du|\leq c\left(\miint{Q_r(z_0)}\left[H(t,|Du|)+1\right]dz\right)^{\frac12}
    \end{equation}
    for every $Q_r(z_0)\Subset\Omega_T$.
\end{theorem}

For $N=1$ and $q<p+\frac4n$, a bound for $Du$ in terms of the $L^q$-norm of $Du$ follows from \cite[Theorem 1.2]{BDM2013}, with a constant and an exponent which blow up as $q\to p+\frac4n$. The estimate \eqref{eq : Lipschitz estimate}, in contrast, holds for all $2\leq p<q<\infty$, involves only the double phase energy, and its constant depends only on $n$ and $q$. For $a\equiv0$ it reduces to the gradient bound known for the degenerate $p$-Laplace system, see \cite[Chapter VIII]{1993_Degenerate_parabolic_equations_DiBenedetto}.

For the H\"older continuity of the gradient we assume that the coefficient is H\"older continuous in time, that is, there exist $\beta\in(0,1]$ and $[a]_\beta\geq0$ such that
\begin{equation}\label{cond : Holder in time}
    |a(t_1)-a(t_2)|\leq [a]_\beta|t_1-t_2|^\beta\qquad\text{for all } t_1,t_2\in(0,T).
\end{equation}
Such a coefficient is bounded. For $\gamma\in(0,1)$ and a cylinder $Q\subset\Omega_T$, we say that $f\in C^{\gamma,\frac{\gamma}{2}}(Q)$ if $f$ is bounded in $Q$ and
$$
[f]_{\gamma;Q}\coloneq\sup_{\substack{(x_1,t_1),(x_2,t_2)\in Q\\ (x_1,t_1)\neq(x_2,t_2)}}\frac{|f(x_1,t_1)-f(x_2,t_2)|}{\max\left\{|x_1-x_2|^\gamma,|t_1-t_2|^{\frac\gamma2}\right\}}<\infty.
$$
Our main result is the following.
\begin{theorem}\label{thm : Holder}
    Let $a:(0,T)\to[0,\infty)$ satisfy \eqref{cond : Holder in time} and let $u$ be a weak solution to \eqref{eq : main equation}. Then $Du$ is locally H\"older continuous in $\Omega_T$. More precisely, $Du\in C^{\gamma,\frac\gamma2}_{\loc}(\Omega_T,\mr^{N\times n})$ with
    $$
    \gamma\coloneq\frac{8\beta\beta_0}{p(n+2+2\beta_0)+8\beta},
    $$
    where $\beta_0=\beta_0(n,N,p,q)\in(0,1)$ is the exponent in Lemma \ref{lem : autonomous}. Moreover, if $Q_{4R}(z_0)\Subset\Omega_T$, then $[Du]_{\gamma;Q_R(z_0)}$ is bounded by a constant depending only on $n,N,p,q,\beta,[a]_\beta,\|a\|_{L^\infty(0,T)},R$ and $\miint{Q_{4R}(z_0)}[H(t,|Du|)+1]\,dz$.
\end{theorem}

Let us comment on these results.
\begin{remark}\label{rem : no gap}
    None of Theorems \ref{thm : second order}, \ref{thm : Lipschitz} and \ref{thm : Holder} involves a condition relating $p$ and $q$. The constants in \eqref{eq : second order estimate} and \eqref{eq : Lipschitz estimate} do not depend on $a$ at all, and Theorems \ref{thm : second order} and \ref{thm : Lipschitz} apply even to discontinuous coefficients such as $a=\chi_{(t_1,T)}$ with $t_1\in(0,T)$, for which the system changes from the $p$-Laplace system to a $(p,q)$-Laplace system at time $t_1$.
\end{remark}

\begin{remark}\label{rem : comparison argument}
    Theorem \ref{thm : Holder} is obtained by combining Theorem \ref{thm : Lipschitz} with a comparison argument. Once the gradient is bounded, the comparison argument does not require a gap condition either, since the oscillation of the coefficient can then be made small by taking small cylinders. The exponent $\gamma$ in Theorem \ref{thm : Holder} comes from an interpolation argument, and we do not expect it to be optimal.
\end{remark}

Let us now describe the ideas of the proofs. The main difficulty is that the coefficient is merely measurable in time, so that the system can neither be differentiated in time nor compared with a frozen system at the level of Theorems \ref{thm : second order} and \ref{thm : Lipschitz}. We therefore work with the system differentiated in the spatial variables, where the two observations above come into play. Theorem \ref{thm : second order} is proved by the difference quotient technique. By the translation invariance, both terms of the difference $\mathcal{A}(t,Du(x+he_k,t))-\mathcal{A}(t,Du(x,t))$, where $e_k$ is a coordinate unit vector, carry the same coefficient $a(t)$, which can be taken out of the spatial integrals, and the energy estimate closes on each time slice. For coefficients depending on $x$, a shifted weight $a(x,t)|Du(x+he_k,t)|^q$ appears at this point, see Remark \ref{rem : x-dependent}. Theorem \ref{thm : Lipschitz} is proved by a Moser iteration for the gradient in the spirit of DiBenedetto and Friedman \cite{DiBenedettoFriedman1984,1993_Degenerate_parabolic_equations_DiBenedetto} and of Uhlenbeck \cite{Uhlenbeck1977}. After differentiation in $x$, no derivative of the coefficient appears, and both phases are treated simultaneously through the primitives of $\sigma\mapsto\mathcal{G}(t,\sigma)^{\frac12}\sigma^m$ with $m\geq0$, to which the Sobolev inequality is applied on each time slice. The resulting estimates depend only on \eqref{cond : Lieberman condition} and on the normalization $\mathcal{G}(t,s)\geq1$ for $s\geq1$. This is possible because the ellipticity and the growth of the differentiated system are both controlled by $\mathcal{G}(t,|Du|)$, whereas for general vector fields with $(p,q)$-growth they are of different orders, which leads to the bound $q<p+\frac4n$ in \cite{BDM2013}. The computations are justified by approximation with regularized systems, in which the growth is truncated at large gradients and the second phase carries the coefficient $a(t)+\varepsilon$, so that no continuity of $a$ is needed.

Finally, Theorem \ref{thm : Holder} is proved by comparing $u$ with solutions of the systems obtained by freezing the coefficient at a fixed time. After normalization at the level of the gradient bound, the frozen systems belong to a family of autonomous systems with Orlicz growth whose gradient regularity is uniform in the phase parameter by the work of Ok, Scilla and Stroffolini \cite{OkScillaStroffolini2024}, and the comparison error is controlled by the oscillation of $a$ in time.

The rest of the paper is organized as follows. In Section \ref{section 2}, we fix notation and collect algebraic inequalities and facts on Steklov averages. Section \ref{section 3} contains the proof of Theorem \ref{thm : second order} and the higher integrability of the gradient. In Section \ref{section 4}, we prove an a priori Lipschitz estimate for systems with Uhlenbeck structure, complete the proof of Theorem \ref{thm : Lipschitz} by approximation, and show that the time derivative of a weak solution is locally square integrable. Section \ref{section 5} is devoted to the proof of Theorem \ref{thm : Holder}.

\section{\bf Preliminaries}\label{section 2}
For $\xi=(\xi_i^\alpha),\eta=(\eta_i^\alpha)\in\mr^{N\times n}$ we write $\xi\cdot\eta\coloneq\sum_{\alpha,i}\xi_i^\alpha\eta_i^\alpha$ and $|\xi|^2\coloneq\xi\cdot\xi$, and for $v\in\mr^N$ and $w\in\mr^n$ we denote by $v\otimes w\in\mr^{N\times n}$ the matrix with entries $v^\alpha w_i$. For $u:\Omega_T\to\mr^N$, $Du$ denotes the spatial gradient, and for $F:\Omega_T\to\mr^{N\times n}$, $\operatorname{div}F$ denotes the vector with components $\sum_iD_iF_i^\alpha$. For a measurable set $E\subset\mr^{n+1}$ with $0<|E|<\infty$ and $f\in L^1(E)$, we write
$$
f_E\coloneq\miint{E}f\,dz\coloneq\frac1{|E|}\iints{E}f\,dz,
$$
where $|E|$ is the $(n+1)$-dimensional Lebesgue measure of $E$. The parabolic boundary of a cylinder $B\times(t_1,t_2)$ is $(\partial B\times(t_1,t_2))\cup(\overline{B}\times\{t_1\})$. For measurable functions, $\sup$ and $\osc$ denote the essential supremum and the essential oscillation. Throughout the paper, $c$ and $C$ denote positive constants which may change from line to line, and their dependence on the parameters is indicated when it matters.

We use the following standard algebraic inequalities, see for instance \cite{Hamburger1992}. For $m\geq2$ and $\mu\geq0$, there exist constants $0<c\leq1\leq C$ depending only on $m$ such that
\begin{align}
\left((\mu^2+|\xi|^2)^{\frac{m-2}{2}}\xi-(\mu^2+|\eta|^2)^{\frac{m-2}{2}}\eta\right)\cdot(\xi-\eta)&\geq c\,(\mu^2+|\xi|^2+|\eta|^2)^{\frac{m-2}{2}}|\xi-\eta|^2, \label{eq : monotonicity}\\
\left|(\mu^2+|\xi|^2)^{\frac{m-2}{2}}\xi-(\mu^2+|\eta|^2)^{\frac{m-2}{2}}\eta\right|&\leq C\,(\mu^2+|\xi|^2+|\eta|^2)^{\frac{m-2}{2}}|\xi-\eta|, \label{eq : growth}\\
c\,(|\xi|+|\eta|)^{m-2}|\xi-\eta|^2\leq|V_m(\xi)-V_m(\eta)|^2&\leq C\,(|\xi|+|\eta|)^{m-2}|\xi-\eta|^2 \label{eq : V equivalence}
\end{align}
for all $\xi,\eta\in\mr^{N\times n}$. Combining \eqref{eq : monotonicity}, \eqref{eq : growth} and \eqref{eq : V equivalence} with $\mu=0$, and using $|\xi|+|\eta|\geq|\xi-\eta|$, we obtain
\begin{equation}\label{eq : V consequences}
\begin{aligned}
&\left(|\xi|^{m-2}\xi-|\eta|^{m-2}\eta\right)\cdot(\xi-\eta)\geq c\,|V_m(\xi)-V_m(\eta)|^2\geq c\,|\xi-\eta|^m,\\
&\left||\xi|^{m-2}\xi-|\eta|^{m-2}\eta\right|\leq C\,(|\xi|+|\eta|)^{\frac{m-2}{2}}|V_m(\xi)-V_m(\eta)|
\end{aligned}
\end{equation}
with constants depending only on $m$.

We next record the properties of vector fields with Uhlenbeck structure. Let $\mathcal{G}:[0,\infty)\to[0,\infty)$ be a continuously differentiable function with
\begin{equation}\label{cond : Uhlenbeck}
0\leq s\,\mathcal{G}'(s)\leq(q-2)\,\mathcal{G}(s)\qquad\text{for all } s\geq0,
\end{equation}
and set $\mathcal{A}_{\mathcal{G}}(\xi)\coloneq\mathcal{G}(|\xi|)\xi$. Then $\mathcal{A}_{\mathcal{G}}\in C^1(\mr^{N\times n},\mr^{N\times n})$, and for $\xi\neq0$ and $\zeta\in\mr^{N\times n}$,
$$
\partial_\xi\mathcal{A}_{\mathcal{G}}(\xi)\zeta=\mathcal{G}(|\xi|)\zeta+\frac{\mathcal{G}'(|\xi|)}{|\xi|}(\xi\cdot\zeta)\,\xi,
$$
while $\partial_\xi\mathcal{A}_{\mathcal{G}}(0)\zeta=\mathcal{G}(0)\zeta$. Hence $\partial_\xi\mathcal{A}_{\mathcal{G}}(\xi)$ is symmetric with eigenvalues $\mathcal{G}(|\xi|)$ and $\mathcal{G}(|\xi|)+|\xi|\mathcal{G}'(|\xi|)$, and \eqref{cond : Uhlenbeck} gives
\begin{equation}\label{eq : eigenvalues}
\mathcal{G}(|\xi|)|\zeta|^2\leq\partial_\xi\mathcal{A}_{\mathcal{G}}(\xi)\zeta\cdot\zeta\leq(q-1)\mathcal{G}(|\xi|)|\zeta|^2,
\qquad
|\partial_\xi\mathcal{A}_{\mathcal{G}}(\xi)\zeta|\leq(q-1)\mathcal{G}(|\xi|)|\zeta|.
\end{equation}
Integrating along the segment from $\eta$ to $\xi$, we obtain
\begin{equation}\label{eq : Uhlenbeck monotonicity}
\left(\mathcal{A}_{\mathcal{G}}(\xi)-\mathcal{A}_{\mathcal{G}}(\eta)\right)\cdot(\xi-\eta)\geq\int_0^1\mathcal{G}(|\eta+\sigma(\xi-\eta)|)\,d\sigma\,|\xi-\eta|^2,
\end{equation}
together with the analogous upper bound for $|\mathcal{A}_{\mathcal{G}}(\xi)-\mathcal{A}_{\mathcal{G}}(\eta)|$ with the additional factor $q-1$. We note that \eqref{cond : Uhlenbeck} also implies
\begin{equation}\label{eq : almost increasing}
\left(\frac{\sigma}{s}\right)^{q-2}\mathcal{G}(s)\leq\mathcal{G}(\sigma)\leq\mathcal{G}(s)\qquad\text{for } 0<\sigma\leq s,
\end{equation}
provided $\mathcal{G}>0$, since $\log\mathcal{G}$ is nondecreasing and $s\mapsto\log\mathcal{G}(s)-(q-2)\log s$ is nonincreasing.

Finally, we recall the Steklov formulation of \eqref{eq : main equation}. For $f\in L^1(\Omega_T)$ and $0<h<T$, the Steklov average of $f$ is
$$
[f]_h(x,t)\coloneq\frac1h\int_t^{t+h}f(x,s)\,ds\qquad\text{for } t\in(0,T-h).
$$
Let $u$ be a weak solution to \eqref{eq : main equation}. Since $u\in L^q(0,T;W^{1,q}(\Omega,\mr^N))$ and $|\mathcal{A}(t,\xi)|\leq|\xi|^{p-1}+\|a\|_{L^\infty(0,T)}|\xi|^{q-1}$, the flux $\mathcal{A}(\cdot,Du)$ belongs to $L^{q'}(\Omega_T)$ with $q'\coloneq\frac{q}{q-1}$. A standard argument, see \cite[Chapter II]{1993_Degenerate_parabolic_equations_DiBenedetto}, shows that
\begin{equation}\label{eq : Steklov}
\int_\Omega\left[\partial_t[u]_h(x,t)\cdot\varphi(x)+[\mathcal{A}(\cdot,Du)]_h(x,t)\cdot D\varphi(x)\right]dx=0
\end{equation}
for every $\varphi\in W^{1,q}_0(\Omega,\mr^N)$ and almost every $t\in(0,T-h)$. Since $W^{1,q}_0(\Omega,\mr^N)$ is separable, the exceptional set of times can be chosen independently of $\varphi$, so that \eqref{eq : Steklov} may be integrated in time against test functions which depend measurably on $t$. Moreover, for every $[t_1,t_2]\subset(0,T)$, we have
$$
\begin{aligned}
&[u]_h\to u\ \ \text{in } C([t_1,t_2];L^2(\Omega,\mr^N)),\qquad [Du]_h\to Du\ \ \text{in } L^q(\Omega\times(t_1,t_2)),\\
&[\mathcal{A}(\cdot,Du)]_h\to\mathcal{A}(\cdot,Du)\ \ \text{in } L^{q'}(\Omega\times(t_1,t_2))
\end{aligned}
$$
as $h\to0$. Note also that the weak formulation gives $u_t\in L^{q'}(0,T;W^{-1,q'}(\Omega,\mr^N))$.

\section{\bf Second order estimates and higher integrability}\label{section 3}
In this section we prove Theorem \ref{thm : second order} and derive the higher integrability of the gradient. The proof uses only the invariance of \eqref{eq : main equation} under spatial translations, which holds because the coefficient does not depend on $x$.

\begin{proof}[Proof of Theorem \ref{thm : second order}]
Fix $Q_{2r}(z_0)\Subset\Omega_T$, a coordinate unit vector $e_k$ with $k\in\{1,\dots,n\}$, and $h\in\mr$ with $0<|h|<r/4$. For a function $f$ defined on $\Omega_T$, we write
$$
\Delta_hf(x,t)\coloneq f(x+he_k,t)-f(x,t)
$$
whenever $x,x+he_k\in\Omega$. We choose $\eta\in C_0^\infty(B_{3r/2}(x_0))$ with $0\leq\eta\leq1$, $\eta\equiv1$ in $B_r(x_0)$ and $|D\eta|\leq 4/r$, and we let $\zeta:\mr\to[0,1]$ be the Lipschitz function with $\zeta(t)=0$ for $t\leq t_0-4r^2$ and $\zeta(t)=1$ for $t\geq t_0-r^2$ which is linear in between, so that $0\leq\zeta'\leq 1/(3r^2)$. For small $\delta>0$, let $[u]_\delta$ be the Steklov average of $u$ with parameter $\delta$. For almost every $t\in\Lambda_{2r}(t_0)$, the function
$$
\varphi(\cdot,t)\coloneq\zeta(t)\,\Delta_{-h}\left(\eta^2\Delta_h[u]_\delta(\cdot,t)\right)
$$
belongs to $W^{1,q}_0(B_{2r}(x_0),\mr^N)$, since $\eta^2\Delta_h[u]_\delta$ is supported in $B_{3r/2}(x_0)$ and $|h|<r/4$. We use it as a test function in \eqref{eq : Steklov} and integrate over $(t_0-4r^2,\tau)$ with $\tau\in\Lambda_{2r}(t_0)$. The discrete integration by parts formula
$$
\int f\cdot\Delta_{-h}g\,dx=\int\Delta_hf\cdot g\,dx,
$$
valid for $g$ supported in $B_{7r/4}(x_0)$, together with the fact that $\Delta_h$ commutes with $D$ and with Steklov averages, gives
$$
\int_{t_0-4r^2}^{\tau}\int_{B_{2r}(x_0)}\zeta\eta^2\,\partial_t\Delta_h[u]_\delta\cdot\Delta_h[u]_\delta\,dx\,dt
+\iints{Q^\tau}\zeta\,[\Delta_h\mathcal{A}]_\delta\cdot D\left(\eta^2\Delta_h[u]_\delta\right)dz=0,
$$
where $Q^\tau\coloneq B_{2r}(x_0)\times(t_0-4r^2,\tau)$ and
$$
\Delta_h\mathcal{A}(x,t)\coloneq\mathcal{A}(t,Du(x+he_k,t))-\mathcal{A}(t,Du(x,t)).
$$
We stress that both terms of $\Delta_h\mathcal{A}$ are evaluated at the same time $t$, and therefore they carry the same coefficient $a(t)$. Since $\zeta(t_0-4r^2)=0$, the first integral equals
$$
\frac12\int_{B_{2r}(x_0)}\zeta\eta^2|\Delta_h[u]_\delta|^2(x,\tau)\,dx-\frac12\iints{Q^\tau}\zeta'\eta^2|\Delta_h[u]_\delta|^2\,dz.
$$
Letting $\delta\to0$ and using the convergence properties of Steklov averages stated after \eqref{eq : Steklov}, we arrive at
\begin{equation}\label{eq : difference identity}
\begin{aligned}
&\frac12\int_{B_{2r}(x_0)}\zeta\eta^2|\Delta_hu|^2(x,\tau)\,dx+\iints{Q^\tau}\zeta\eta^2\,\Delta_h\mathcal{A}\cdot\Delta_hDu\,dz\\
&\qquad=\frac12\iints{Q^\tau}\zeta'\eta^2|\Delta_hu|^2\,dz-\iints{Q^\tau}\zeta\,\Delta_h\mathcal{A}\cdot\left(\Delta_hu\otimes D\eta^2\right)dz
\end{aligned}
\end{equation}
for every $\tau\in\Lambda_{2r}(t_0)$.

We now estimate the terms in \eqref{eq : difference identity}. At a point $(x,t)$, we write $\xi_0\coloneq Du(x,t)$, $\xi_1\coloneq Du(x+he_k,t)$ and $S\coloneq|\xi_0|+|\xi_1|$. Applying \eqref{eq : V consequences} with $m=p$ and $m=q$, and using $a(t)\geq0$, we obtain
\begin{equation}\label{eq : difference monotone}
\Delta_h\mathcal{A}\cdot\Delta_hDu\geq c_0\left(|\Delta_hV_p(Du)|^2+a(t)|\Delta_hV_q(Du)|^2\right)
\end{equation}
and
$$
|\Delta_h\mathcal{A}|\leq C\left(S^{\frac{p-2}{2}}|\Delta_hV_p(Du)|+a(t)S^{\frac{q-2}{2}}|\Delta_hV_q(Du)|\right),
$$
where $c_0$ and $C$ depend only on $p$ and $q$. By Young's inequality,
$$
\begin{aligned}
\left|\zeta\,\Delta_h\mathcal{A}\cdot\left(\Delta_hu\otimes D\eta^2\right)\right|
&\leq\frac{c_0}2\zeta\eta^2\left(|\Delta_hV_p(Du)|^2+a(t)|\Delta_hV_q(Du)|^2\right)\\
&\quad+C\zeta|D\eta|^2|\Delta_hu|^2\left(S^{p-2}+a(t)S^{q-2}\right).
\end{aligned}
$$
The last term is controlled on each time slice. For almost every $t$ and $m\in\{p,q\}$, we have
$$
\int_{B_{3r/2}(x_0)}S^m\,dx\leq2^m\int_{B_{7r/4}(x_0)}|Du(x,t)|^m\,dx
$$
and
$$
\int_{B_{3r/2}(x_0)}|\Delta_hu(x,t)|^m\,dx\leq|h|^m\int_{B_{7r/4}(x_0)}|Du(x,t)|^m\,dx.
$$
If $m>2$, these two inequalities and H\"older's inequality with exponents $\frac{m}{m-2}$ and $\frac m2$ give
\begin{equation}\label{eq : slice estimate}
\int_{B_{3r/2}(x_0)}S^{m-2}|\Delta_hu|^2\,dx\leq C|h|^2\int_{B_{7r/4}(x_0)}|Du(x,t)|^m\,dx,
\end{equation}
and \eqref{eq : slice estimate} is immediate for $m=2$. Since $a$ does not depend on $x$, we may multiply \eqref{eq : slice estimate} with $m=q$ by $a(t)$ and obtain
$$
\int_{B_{3r/2}(x_0)}a(t)S^{q-2}|\Delta_hu|^2\,dx\leq C|h|^2\int_{B_{7r/4}(x_0)}a(t)|Du(x,t)|^q\,dx.
$$
This is the only place where the structure of the coefficient enters. For the first term on the right-hand side of \eqref{eq : difference identity}, we use $\zeta'\leq 1/(3r^2)$ and
$$
\int_{B_{3r/2}(x_0)}|\Delta_hu|^2\,dx\leq|h|^2\int_{B_{7r/4}(x_0)}|Du|^2\,dx.
$$
Inserting these estimates into \eqref{eq : difference identity}, absorbing the first term of the preceding Young inequality into the left-hand side and recalling that $\zeta\eta^2=1$ in $Q_r(z_0)$, we find
$$
\begin{aligned}
&\int_{B_r(x_0)}|\Delta_hu|^2(x,\tau)\,dx+\int_{t_0-r^2}^{\tau}\int_{B_r(x_0)}\left(|\Delta_hV_p(Du)|^2+a(t)|\Delta_hV_q(Du)|^2\right)dx\,dt\\
&\qquad\leq\frac{C|h|^2}{r^2}\iints{Q_{2r}(z_0)}\left[|Du|^2+H(t,|Du|)\right]dz
\end{aligned}
$$
for every $\tau\in\Lambda_r(t_0)$, where $C=C(p,q)$.

Consequently the difference quotients of $V_p(Du)$ and of $a^{\frac12}V_q(Du)$ in the direction $e_k$ are bounded in $L^2(Q_r(z_0))$ uniformly in $h$, and for every $\tau\in\Lambda_r(t_0)$ the difference quotients $\Delta_hu(\cdot,\tau)/h$ are bounded in $L^2(B_r(x_0))$. By the difference quotient characterization of Sobolev functions, the weak derivatives $D_kV_p(Du)$ and $D_k(a^{\frac12}V_q(Du))$ exist in $L^2(Q_r(z_0))$, and $D_ku(\cdot,\tau)\in L^2(B_r(x_0))$ for every $\tau\in\Lambda_r(t_0)$, with the same bounds. By Fubini's theorem, for almost every $t\in\Lambda_r(t_0)$ the functions $V_p(Du(\cdot,t))$ and $a(t)^{\frac12}V_q(Du(\cdot,t))$ then belong to $W^{1,2}(B_r(x_0),\mr^{N\times n})$, and their spatial gradients are the restrictions of the space-time derivatives. Summing over $k=1,\dots,n$ proves \eqref{eq : second order estimate}.
\end{proof}

We point out where the preceding argument breaks down for coefficients depending on $x$.
\begin{remark}\label{rem : x-dependent}
    If $a=a(x,t)$, the difference $\Delta_h\mathcal{A}$ contains the additional term $(a(x+he_k,t)-a(x,t))|\xi_1|^{q-2}\xi_1$, which has no sign, and the slice estimate \eqref{eq : slice estimate} with $m=q$ would have to be carried out with the weight $a(x,t)$ against $|Du(x+he_k,t)|^q$. Replacing this shifted weight by $a(x+he_k,t)$ requires control of the oscillation of $a$, and this is where gap conditions enter the difference quotient technique.
\end{remark}

As a consequence of Theorem \ref{thm : second order}, the gradient enjoys higher integrability with an explicit exponent.
\begin{corollary}\label{cor : higher integrability}
    Let $a\in L^\infty(0,T)$ be nonnegative and let $u$ be a weak solution to \eqref{eq : main equation}. Then
    \begin{equation}\label{eq : explicit higher integrability}
        \miint{Q_r(z_0)}\left[|Du|^{p+\frac4n}+a(t)|Du|^{q+\frac4n}\right]dz\leq c\left(\miint{Q_{2r}(z_0)}\left[|Du|^2+H(t,|Du|)\right]dz\right)^{1+\frac2n}
    \end{equation}
    for every $Q_{2r}(z_0)\Subset\Omega_T$, where $c=c(n,p,q)$.
\end{corollary}
\begin{proof}
Assume first that $n\geq3$, and let $2^*\coloneq\frac{2n}{n-2}$. For almost every $t\in\Lambda_r(t_0)$, we have $V_p(Du(\cdot,t))\in W^{1,2}(B_r(x_0))$ by Theorem \ref{thm : second order}, and the Sobolev inequality on $B_r=B_r(x_0)$ reads
$$
\|f\|_{L^{2^*}(B_r)}^2\leq c(n)\left(\|Df\|_{L^2(B_r)}^2+r^{-2}\|f\|_{L^2(B_r)}^2\right)\qquad\text{for } f\in W^{1,2}(B_r).
$$
H\"older's inequality with exponents $\frac{n}{n-2}$ and $\frac n2$ yields
$$
\begin{aligned}
\int_{B_r}|Du|^{p+\frac4n}\,dx&=\int_{B_r}|V_p(Du)|^2|Du|^{\frac4n}\,dx\\
&\leq c\left(\int_{B_r}|DV_p(Du)|^2\,dx+\frac1{r^2}\int_{B_r}|Du|^p\,dx\right)\left(\int_{B_r}|Du|^2\,dx\right)^{\frac2n}.
\end{aligned}
$$
The same computation applies to $a(t)|Du|^{q+\frac4n}=|a(t)^{\frac12}V_q(Du)|^2|Du|^{\frac4n}$, where the Sobolev inequality is applied to $a(t)^{\frac12}V_q(Du(\cdot,t))$. This is possible because $a(t)$ does not depend on $x$. We obtain
$$
\begin{aligned}
&\int_{B_r}a(t)|Du|^{q+\frac4n}\,dx\\
&\qquad\leq c\left(\int_{B_r}a(t)|DV_q(Du)|^2\,dx+\frac1{r^2}\int_{B_r}a(t)|Du|^q\,dx\right)\left(\int_{B_r}|Du|^2\,dx\right)^{\frac2n}.
\end{aligned}
$$
We integrate these two inequalities over $\Lambda_r(t_0)$ and apply \eqref{eq : second order estimate} to both factors. This gives
$$
\begin{aligned}
&\iints{Q_r(z_0)}\left[|Du|^{p+\frac4n}+a(t)|Du|^{q+\frac4n}\right]dz\\
&\qquad\leq c\left(\frac1{r^2}\iints{Q_{2r}(z_0)}\left[|Du|^2+H(t,|Du|)\right]dz\right)^{1+\frac2n}.
\end{aligned}
$$
Dividing by $|Q_r(z_0)|$ and using $|Q_{2r}(z_0)|^{1+\frac2n}r^{-2-\frac4n}=c(n)|Q_r(z_0)|$, we obtain \eqref{eq : explicit higher integrability}.

If $n=2$, we replace the Sobolev inequality in $W^{1,2}$ by the one in $W^{1,1}$. Let $m\in\{p,q\}$ and let $t\in\Lambda_r(t_0)$ be such that $V_m(Du(\cdot,t))\in W^{1,2}(B_r)$, where for $m=q$ we consider only the times with $a(t)>0$. We write $W\coloneq|Du(\cdot,t)|$ and $F\coloneq|V_m(Du(\cdot,t))|^{\frac{m+2}{m}}=W^{\frac{m+2}{2}}$. By the chain rule, which may be justified by truncation, $|DF|\leq\frac{m+2}{m}W|DV_m(Du)|$, and the Cauchy--Schwarz inequality gives
$$
\int_{B_r}F\,dx\leq\left(\int_{B_r}W^2\,dx\right)^{\frac12}\left(\int_{B_r}W^m\,dx\right)^{\frac12}
$$
and
$$
\int_{B_r}|DF|\,dx\leq\frac{m+2}{m}\left(\int_{B_r}W^2\,dx\right)^{\frac12}\left(\int_{B_r}|DV_m(Du)|^2\,dx\right)^{\frac12}.
$$
We combine these estimates with the Sobolev inequality
$$
\|F\|_{L^2(B_r)}\leq c\left(\|DF\|_{L^1(B_r)}+r^{-1}\|F\|_{L^1(B_r)}\right)
$$
in two dimensions and obtain
$$
\begin{aligned}
\int_{B_r}|Du|^{m+2}\,dx&=\|F\|^2_{L^2(B_r)}\\
&\leq c\left(\int_{B_r}|DV_m(Du)|^2\,dx+\frac1{r^2}\int_{B_r}|Du|^m\,dx\right)\int_{B_r}|Du|^2\,dx.
\end{aligned}
$$
For $m=p$ and, after multiplication by $a(t)$, for $m=q$, these are the two estimates above with $n=2$, and the rest of the proof is the same.
\end{proof}

In particular, the function $H(t,|Du|)^{1+\frac{4}{nq}}$ belongs to $L^1_{\loc}(\Omega_T)$, because $p(1+\frac4{nq})\leq p+\frac4n$ and $a(t)^{1+\frac4{nq}}\leq\|a\|_{L^\infty(0,T)}^{\frac{4}{nq}}a(t)$. For parabolic problems with $p$-growth, the gain from $p$ to $p+\frac4n$ is the natural one, and it corresponds to the upper bound $q<p+\frac4n$ in \cite{BDM2013}, see \cite[Remark 1.4]{BDM2013}. Gehring type arguments yield instead a small and implicit exponent, see \cite{Kinnunen2000} for the $p$-Laplace system and \cite{2023_Gradient_Higher_Integrability_for_Degenerate_Parabolic_Double-Phase_Systems} for double phase systems under a gap bound.

\section{\bf Lipschitz bounds}\label{section 4}
In this section we prove Theorem \ref{thm : Lipschitz} and deduce the local square integrability of the time derivative. We first prove an a priori estimate for systems with Uhlenbeck structure and bounded ellipticity, whose constant depends only on $n$ and $q$, and then we remove the additional assumptions by approximation.

Let $\mathcal{G}:(0,T)\times[0,\infty)\to[0,\infty)$ be a function such that $t\mapsto\mathcal{G}(t,s)$ is measurable for every $s\geq0$, $s\mapsto\mathcal{G}(t,s)$ belongs to $C^1([0,\infty))$ for every $t\in(0,T)$, and
\begin{equation}\label{cond : G}
0\leq s\,\partial_s\mathcal{G}(t,s)\leq(q-2)\,\mathcal{G}(t,s)\quad\text{for } s\geq0,
\qquad
\mathcal{G}(t,s)\geq1\quad\text{for } s\geq1.
\end{equation}
We also assume that there exist constants $0<\lambda_0\leq\lambda_1<\infty$ such that
\begin{equation}\label{cond : bounded ellipticity}
\lambda_0\leq\mathcal{G}(t,s)\leq\lambda_1\qquad\text{for all } t\in(0,T) \text{ and } s\geq0,
\end{equation}
and we set $\mathcal{A}_{\mathcal{G}}(t,\xi)\coloneq\mathcal{G}(t,|\xi|)\xi$. The model field \eqref{eq : vector field} corresponds to $\mathcal{G}(t,s)=s^{p-2}+a(t)s^{q-2}$, which satisfies \eqref{cond : G} for $s>0$ but in general is neither differentiable at $s=0$ nor bounded as in \eqref{cond : bounded ellipticity}. Assumption \eqref{cond : bounded ellipticity} is used only qualitatively, and the constants below depend neither on $\lambda_0$ nor on $\lambda_1$.

The following a priori estimate is the core of the proof of Theorem \ref{thm : Lipschitz}.
\begin{proposition}\label{prop : a priori}
    Let $\mathcal{G}$ satisfy \eqref{cond : G} and \eqref{cond : bounded ellipticity}, let $Q=B\times I$ be a cylinder, and let
    $$
    v\in C(\overline{I};L^2(B,\mr^N))\cap L^2(I;W^{1,2}(B,\mr^N))
    $$
    be a weak solution to $v_t-\operatorname{div}\mathcal{A}_{\mathcal{G}}(t,Dv)=0$ in $Q$, that is, the identity in Definition \ref{def : weak solution} holds with $\mathcal{A}_{\mathcal{G}}$ in place of $\mathcal{A}$ and $Q$ in place of $\Omega_T$. Then $Dv\in L^\infty_{\loc}(Q,\mr^{N\times n})$, and there exists a constant $c=c(n,q)\geq1$ such that
    \begin{equation}\label{eq : a priori estimate}
        \sup_{Q_{\rho/2}(z_0)}|Dv|\leq c\left(\miint{Q_\rho(z_0)}\left[\mathcal{G}(t,|Dv|)|Dv|^2+1\right]dz\right)^{\frac12}
    \end{equation}
    for every $Q_\rho(z_0)\Subset Q$.
\end{proposition}
\begin{proof}
We first show that $v$ has second spatial derivatives and derive the differentiated system. By \eqref{eq : eigenvalues}, \eqref{eq : Uhlenbeck monotonicity} and \eqref{cond : bounded ellipticity}, we have
$$
\left(\mathcal{A}_{\mathcal{G}}(t,\xi)-\mathcal{A}_{\mathcal{G}}(t,\eta)\right)\cdot(\xi-\eta)\geq\lambda_0|\xi-\eta|^2,
\qquad
|\mathcal{A}_{\mathcal{G}}(t,\xi)-\mathcal{A}_{\mathcal{G}}(t,\eta)|\leq(q-1)\lambda_1|\xi-\eta|
$$
for all $t$ and $\xi,\eta\in\mr^{N\times n}$. In particular, $\mathcal{A}_{\mathcal{G}}(\cdot,Dv)\in L^2(Q)$, and the Steklov formulation \eqref{eq : Steklov} holds for $v$ with test functions in $W^{1,2}_0(B,\mr^N)$. Repeating the difference quotient argument in the proof of Theorem \ref{thm : second order}, with the two preceding inequalities in place of \eqref{eq : V consequences}, we find
$$
\iints{Q'}|\Delta_hDv|^2\,dz\leq C|h|^2
$$
for every $Q'\Subset Q$, with $C$ depending on $\lambda_0$, $\lambda_1$, $q$, $Q'$ and $\|Dv\|_{L^2(Q)}$. Thus $D_kv\in L^2_{\loc}(I;W^{1,2}_{\loc}(B,\mr^N))$ for $k=1,\dots,n$. Since $\mathcal{A}_{\mathcal{G}}(t,\cdot)$ is continuously differentiable with bounded derivative, the chain rule gives $D_k[\mathcal{A}_{\mathcal{G}}(t,Dv)]=\mathcal{B}\,DD_kv$ with
$$
\mathcal{B}(z)\coloneq\partial_\xi\mathcal{A}_{\mathcal{G}}(t,Dv(z)),
$$
a measurable field of symmetric linear maps on $\mr^{N\times n}$ which satisfies \eqref{eq : eigenvalues} with $\xi=Dv(z)$ and is bounded by $(q-1)\lambda_1$. Testing the system with $D_k\varphi$ for $\varphi\in C_0^\infty(Q,\mr^N)$ and integrating by parts in $x_k$, we see that $w_k\coloneq D_kv$ satisfies
\begin{equation}\label{eq : differentiated system}
\iints{Q}\left[-w_k\cdot\varphi_t+\mathcal{B}Dw_k\cdot D\varphi\right]dz=0
\end{equation}
for every $\varphi\in C_0^\infty(Q,\mr^N)$. Hence $\partial_tw_k$ belongs locally to $L^2(I;W^{-1,2}(B,\mr^N))$, and the identity \eqref{eq : differentiated system} extends to all test functions $\varphi$ in $L^2(I;W^{1,2}(B,\mr^N))$ which vanish outside a compact subset of $Q$, provided that the first term is replaced by the duality pairing
$$
\int_I\langle\partial_tw_k,\varphi\rangle\,dt.
$$

Next, we derive the energy inequality for powers of $W\coloneq|Dv|$. Fix $Q_\rho(z_0)\Subset Q$ and a function $\eta\in C^\infty(\mr^{n+1})$ with $0\leq\eta\leq1$ such that $\eta(\cdot,t)$ is supported in $B_\rho(x_0)$ for every $t$ and $\eta=0$ for $t\leq t_0-\rho^2$. Let $m\geq0$ and $K\geq1$. If $m=0$, we set $\Phi_K\equiv1$. If $m>0$, we choose a nonincreasing $\zeta\in C^\infty([0,\infty))$ with $\zeta=1$ on $[0,1]$ and $\zeta=0$ on $[2,\infty)$, and we set
$$
\Phi_K(s)\coloneq\int_0^sm\sigma^{m-1}\zeta\Bigl(\frac{\sigma}{K}\Bigr)\,d\sigma,
\qquad
\Psi_K(s)\coloneq\int_0^s\Phi_K(\sigma)\,d\sigma.
$$
Then $\Phi_K$ is nondecreasing, $\Phi_K(s)=s^m$ for $s\leq K$, $\Phi_K$ is constant on $[2K,\infty)$, and $\Phi_K(s)\to s^m$ as $K\to\infty$. The map $F_K(\xi)\coloneq\Phi_K(|\xi|^2)\xi$ belongs to $C^1(\mr^{N\times n},\mr^{N\times n})$ and has bounded derivative, because $|\xi|^2\Phi_K'(|\xi|^2)=m|\xi|^{2m}\zeta(|\xi|^2/K)$ tends to zero as $\xi\to0$ and vanishes for $|\xi|^2\geq2K$, and it is the gradient of the convex function $\xi\mapsto\frac12\Psi_K(|\xi|^2)$. We test \eqref{eq : differentiated system} with $\eta^2\Phi_K(W^2)w_k$ on $B_\rho(x_0)\times(t_0-\rho^2,\tau)$ for $\tau\in\Lambda_\rho(t_0)$ and sum over $k$. Since $w_k$ belongs locally to $L^2(I;W^{1,2}(B,\mr^N))$ and $\partial_tw_k$ belongs locally to $L^2(I;W^{-1,2}(B,\mr^N))$, the function $w_k$ has a representative which is continuous in time with values in $L^2_{\loc}(B,\mr^N)$, and a regularization in time, which is justified because $DF_K$ is bounded and continuous, shows that the parabolic part equals
$$
\begin{aligned}
&\sum_{k=1}^n\int_{t_0-\rho^2}^\tau\left\langle\partial_tw_k,\eta^2\Phi_K(W^2)w_k\right\rangle dt\\
&\qquad=\frac12\int_{B_\rho(x_0)}\eta^2\Psi_K(W^2)(x,\tau)\,dx-\iints{Q_\rho^\tau}\eta\eta_t\Psi_K(W^2)\,dz,
\end{aligned}
$$
where $Q_\rho^\tau\coloneq B_\rho(x_0)\times(t_0-\rho^2,\tau)$. For the elliptic part, we note that $\Phi_K(W^2)w_k$ is the $k$-th column of $F_K(Dv)$, whose spatial derivative is $DF_K(Dv)D^2v$. At points where $Dv\neq0$, this gives
$$
\begin{aligned}
D\left(\eta^2\Phi_K(W^2)w_k\right)&=\eta^2\Phi_K(W^2)Dw_k+\eta^2\Phi_K'(W^2)\,w_k\otimes D(W^2)\\
&\quad+\Phi_K(W^2)\,w_k\otimes D(\eta^2).
\end{aligned}
$$
On the set where $Dv=0$, the derivative of $\eta^2\Phi_K(W^2)w_k$ vanishes if $m>0$, since $DF_K(0)=0$, and equals $\eta^2Dw_k$ if $m=0$, so that the following computation is needed only where $Dv\neq0$.
By \eqref{eq : eigenvalues},
$$
\sum_{k=1}^n\mathcal{B}Dw_k\cdot Dw_k\geq\mathcal{G}(t,W)|D^2v|^2,\qquad |D^2v|^2\coloneq\sum_{k=1}^n|Dw_k|^2.
$$
For the second term, since $Dv\cdot Dw_k=\frac12D_k(W^2)$ and $\sum_{k,\alpha}w_k^\alpha D_iw_k^\alpha=\frac12D_i(W^2)$ at points where $Dv\neq0$, we have for every $X\in\mr^n$
$$
\sum_{k=1}^n\mathcal{B}Dw_k\cdot\left(w_k\otimes X\right)
=\frac{\mathcal{G}(t,W)}{2}D(W^2)\cdot X+\frac{\partial_s\mathcal{G}(t,W)}{2W}\sum_{i,j=1}^nX_iM_{ij}D_j(W^2),
$$
where $M_{ij}\coloneq D_iv\cdot D_jv$. The matrix $(M_{ij})$ is positive semidefinite and $\partial_s\mathcal{G}\geq0$. Choosing $X=D(W^2)$, we see that
$$
\sum_{k=1}^n\mathcal{B}Dw_k\cdot\left(w_k\otimes D(W^2)\right)\geq0,
$$
and this also holds at points where $Dv=0$, since $D(W^2)=0$ almost everywhere on this set. As $\Phi_K'\geq0$, the second term is nonnegative. Finally, by \eqref{eq : eigenvalues} and Young's inequality,
$$
\begin{aligned}
\left|\sum_{k=1}^n\mathcal{B}Dw_k\cdot\left(w_k\otimes D(\eta^2)\right)\right|
&\leq2(q-1)\mathcal{G}(t,W)|D^2v|\,W\eta|D\eta|\\
&\leq\frac12\eta^2\mathcal{G}(t,W)|D^2v|^2+2(q-1)^2\mathcal{G}(t,W)W^2|D\eta|^2.
\end{aligned}
$$
Collecting these estimates, we arrive at
$$
\begin{aligned}
&\frac12\int_{B_\rho(x_0)}\eta^2\Psi_K(W^2)(x,\tau)\,dx+\frac12\iints{Q_\rho^\tau}\eta^2\Phi_K(W^2)\mathcal{G}(t,W)|D^2v|^2\,dz\\
&\qquad\leq\iints{Q_\rho(z_0)}\left[2(q-1)^2\Phi_K(W^2)\mathcal{G}(t,W)W^2|D\eta|^2+\Psi_K(W^2)\eta|\eta_t|\right]dz
\end{aligned}
$$
for every $\tau\in\Lambda_\rho(t_0)$. Letting $K\to\infty$ and using monotone convergence, taking the supremum over $\tau$ in the first term on the left-hand side and letting $\tau\to t_0+\rho^2$ in the second one, we obtain
\begin{equation}\label{eq : energy for powers}
\begin{aligned}
&\frac1{m+1}\sup_{\tau\in\Lambda_\rho(t_0)}\int_{B_\rho(x_0)}\eta^2W^{2m+2}(x,\tau)\,dx+\iints{Q_\rho(z_0)}\eta^2\mathcal{G}(t,W)W^{2m}|D^2v|^2\,dz\\
&\qquad\leq c_1\iints{Q_\rho(z_0)}\left[\mathcal{G}(t,W)W^{2m+2}|D\eta|^2+\frac{W^{2m+2}}{m+1}\eta|\eta_t|\right]dz
\end{aligned}
\end{equation}
with $c_1\coloneq8(q-1)^2$. Both sides may be infinite, but the left-hand side is finite whenever the right-hand side is.

We now apply the Sobolev inequality on time slices. By Kato's inequality, $|DW|\leq|D^2v|$ almost everywhere. For the function
$$
\Gamma_m(t,s)\coloneq\int_0^s\mathcal{G}(t,\sigma)^{\frac12}\sigma^m\,d\sigma,
$$
we have
$$
|D[\Gamma_m(t,W)]|^2=\mathcal{G}(t,W)W^{2m}|DW|^2\leq\mathcal{G}(t,W)W^{2m}|D^2v|^2,
$$
where the chain rule is justified by truncating $W$ at a finite level. Here it is essential that $\mathcal{G}$ does not depend on $x$. Moreover, \eqref{eq : almost increasing} yields
\begin{equation}\label{eq : Gamma bounds}
\frac{\mathcal{G}(t,s)^{\frac12}s^{m+1}}{m+\frac q2}\leq\Gamma_m(t,s)\leq\frac{\mathcal{G}(t,s)^{\frac12}s^{m+1}}{m+1}.
\end{equation}
Let
$$
E_m\coloneq\sup_{\tau\in\Lambda_\rho(t_0)}\int_{B_\rho(x_0)}\eta^2W^{2m+2}(x,\tau)\,dx+\iints{Q_\rho(z_0)}|D(\eta\,\Gamma_m(t,W))|^2\,dz.
$$
Using $|D(\eta\Gamma_m)|^2\leq2\eta^2|D\Gamma_m|^2+2\Gamma_m^2|D\eta|^2$, \eqref{eq : Gamma bounds} and \eqref{eq : energy for powers} multiplied by $m+1$, we find
\begin{equation}\label{eq : Em bound}
E_m\leq4c_1(m+1)\iints{Q_\rho(z_0)}\left[\mathcal{G}(t,W)W^{2m+2}|D\eta|^2+W^{2m+2}\eta|\eta_t|\right]dz.
\end{equation}
Assume first that $n\geq3$. For almost every $t$, the function $\eta\Gamma_m(t,W)$ belongs to $W^{1,2}_0(B_\rho(x_0))$, and H\"older's inequality with exponents $\frac{n}{n-2}$ and $\frac n2$ together with the Sobolev inequality gives
$$
\begin{aligned}
&\int_{B_\rho(x_0)}\eta^2\Gamma_m(t,W)^2\left(\eta^2W^{2m+2}\right)^{\frac2n}dx\\
&\qquad\leq c(n)\int_{B_\rho(x_0)}|D(\eta\Gamma_m(t,W))|^2\,dx\left(\int_{B_\rho(x_0)}\eta^2W^{2m+2}\,dx\right)^{\frac2n}.
\end{aligned}
$$
We set $\kappa\coloneq1+\frac2n$. Integrating in time and using the lower bound in \eqref{eq : Gamma bounds}, we obtain
\begin{equation}\label{eq : Sobolev step}
\iints{Q_\rho(z_0)}\eta^{2\kappa}\mathcal{G}(t,W)W^{2(m+1)\kappa}\,dz\leq c(n)\left(m+q\right)^2E_m^{\kappa}.
\end{equation}

Finally, we iterate. For $\frac\rho2\leq\rho'<\rho''\leq\rho$, we choose $\eta$ as above with $\eta(\cdot,t)$ supported in $B_{\rho''}(x_0)$, $\eta=0$ for $t\leq t_0-(\rho'')^2$, $\eta=1$ in $Q_{\rho'}(z_0)$, $|D\eta|\leq2/(\rho''-\rho')$ and $0\leq\eta_t\leq2/((\rho'')^2-(\rho')^2)\leq2/(\rho''-\rho')^2$. Since $W^{2m+2}\leq\mathcal{G}(t,W)W^{2m+2}+1$ by the second condition in \eqref{cond : G}, the quantity
$$
Y_m(\varrho)\coloneq\iints{Q_{\varrho}(z_0)}\left[\mathcal{G}(t,W)W^{2m+2}+1\right]dz
$$
satisfies, by \eqref{eq : Em bound} and \eqref{eq : Sobolev step} applied with $\rho''$ in place of $\rho$,
$$
\iints{Q_{\rho'}(z_0)}\mathcal{G}(t,W)W^{2(m+1)\kappa}\,dz\leq c(n,q)(m+q)^{2+\kappa}(\rho''-\rho')^{-2\kappa}Y_m(\rho'')^\kappa.
$$
Moreover, $Y_m(\rho'')\geq|Q_{\rho''}(z_0)|$ and $n\kappa=n+2$, so that $|Q_{\rho'}(z_0)|\leq c(n)(\rho''-\rho')^{-2\kappa}Y_m(\rho'')^\kappa$. Therefore
\begin{equation}\label{eq : one step}
Y_{\bar m}(\rho')\leq c_2(m+q)^{2+\kappa}(\rho''-\rho')^{-2\kappa}Y_m(\rho'')^\kappa,\qquad \bar m+1\coloneq\kappa(m+1),
\end{equation}
with $c_2=c_2(n,q)$. Let $\rho_j\coloneq\frac\rho2+2^{-j-1}\rho$, $m_j\coloneq\kappa^j-1$ and $Z_j\coloneq\rho^{-(n+2)}Y_{m_j}(\rho_j)$ for $j\geq0$. Since $Y_0(\rho)<\infty$ by \eqref{cond : bounded ellipticity}, the estimate \eqref{eq : one step} shows inductively that all $Z_j$ are finite. As $(n+2)(\kappa-1)=2\kappa$, $\rho_j-\rho_{j+1}=2^{-j-2}\rho$ and $m_j+q\leq(1+q)\kappa^j$, it follows from \eqref{eq : one step} that
$$
Z_{j+1}\leq c_3b^jZ_j^\kappa,\qquad b\coloneq4^\kappa\kappa^{2+\kappa},
$$
with $c_3=c_3(n,q)\geq1$. By induction, $Z_j\leq(c_3^{s_j}b^{t_j}Z_0)^{\kappa^j}$ with
$$
s_j\coloneq\sum_{i=0}^{j-1}\kappa^{-i-1}\leq\frac1{\kappa-1}\qquad\text{and}\qquad t_j\coloneq\sum_{i=0}^{j-1}i\kappa^{-i-1}\leq\frac1{(\kappa-1)^2},
$$
hence
$$
Z_j\leq(c_4Z_0)^{\kappa^j},\qquad c_4\coloneq c_3^{\frac1{\kappa-1}}b^{\frac1{(\kappa-1)^2}}.
$$
Since $\max\{W,1\}^{2\kappa^j}\leq\mathcal{G}(t,W)W^{2(m_j+1)}+1$ and $Q_{\rho/2}(z_0)\subset Q_{\rho_j}(z_0)$, we conclude that
$$
\left(\miint{Q_{\rho/2}(z_0)}\max\{W,1\}^{2\kappa^j}dz\right)^{\frac1{2\kappa^j}}\leq\left(\frac{\rho^{n+2}Z_j}{|Q_{\rho/2}(z_0)|}\right)^{\frac1{2\kappa^j}}\leq c(n)^{\frac1{2\kappa^j}}(c_4Z_0)^{\frac12}.
$$
Letting $j\to\infty$, we obtain $\sup_{Q_{\rho/2}(z_0)}W\leq(c_4Z_0)^{\frac12}$. This is \eqref{eq : a priori estimate}, because
$$
Z_0=\rho^{-(n+2)}Y_0(\rho)=2|B_1|\miint{Q_\rho(z_0)}\left[\mathcal{G}(t,W)W^2+1\right]dz.
$$

If $n=2$, we regard $v$ as a function of $(x,x_3,t)\in B\times(-2\rho,2\rho)\times I$ which does not depend on $x_3$. This function solves the same system in dimension three, since $\mathcal{A}_{\mathcal{G}}(t,\xi)=\mathcal{G}(t,|\xi|)\xi$, and its gradient has the same norm. Applying the case $n=3$ in the cylinder $B^{(3)}_\rho(x_0,0)\times\Lambda_\rho(t_0)$, where $B^{(3)}_\rho$ is the three-dimensional ball, and using $B^{(3)}_{\rho}(x_0,0)\subset B_\rho(x_0)\times(-\rho,\rho)$ and $|B^{(3)}_\rho|=\frac43\pi\rho^3$ on the right-hand side, we obtain \eqref{eq : a priori estimate} for $n=2$, since by Fubini's theorem the essential supremum of $|Dv|$ over $B^{(3)}_{\rho/2}(x_0,0)\times\Lambda_{\rho/2}(t_0)$ coincides with the one over $Q_{\rho/2}(z_0)$.
\end{proof}

We now turn to the model system. Let $u$ be a weak solution to \eqref{eq : main equation} with a nonnegative coefficient $a\in L^\infty(0,T)$, and fix a cylinder $Q_0\coloneq B_0\times I_0\Subset\Omega_T$ with $I_0=(t_1,t_2)$. For $\varepsilon\in(0,1]$ and $M\geq1$, we define
$$
\mathcal{G}_\varepsilon(t,s)\coloneq(\varepsilon^2+s^2)^{\frac{p-2}{2}}+(a(t)+\varepsilon)(\varepsilon^2+s^2)^{\frac{q-2}{2}},
\qquad
\mathcal{G}_{\varepsilon,M}(t,s)\coloneq\mathcal{G}_\varepsilon(t,\psi_M(s)),
$$
where $\psi_M\in C^\infty([0,\infty))$ is a fixed concave nondecreasing function with $\psi_M(s)=s$ for $0\leq s\leq M$ which is constant on $[3M,\infty)$. We write $\mathcal{A}_\varepsilon\coloneq\mathcal{A}_{\mathcal{G}_\varepsilon}$ and $\mathcal{A}_{\varepsilon,M}\coloneq\mathcal{A}_{\mathcal{G}_{\varepsilon,M}}$. A direct computation shows that $\mathcal{G}_\varepsilon$ satisfies \eqref{cond : G}. Since $\psi_M$ is concave with $\psi_M(0)=0$ and $\psi_M'=1$ on $[0,M]$, we have $s\psi_M'(s)\leq\psi_M(s)\leq\min\{s,3M\}$, and since $\psi_M(s)\geq\min\{s,M\}$, the function $\mathcal{G}_{\varepsilon,M}$ satisfies \eqref{cond : G} as well. Moreover, $\mathcal{G}_{\varepsilon,M}$ satisfies \eqref{cond : bounded ellipticity} with $\lambda_0=\varepsilon^{p-2}$ and $\lambda_1=\sup_t\mathcal{G}_\varepsilon(t,3M)<\infty$. Finally, $\mathcal{G}_{\varepsilon,M}\leq\mathcal{G}_\varepsilon$, since $\mathcal{G}_\varepsilon(t,\cdot)$ is nondecreasing.

In the first approximation step we let the truncation level tend to infinity.
\begin{lemma}\label{lem : truncation limit}
    Let $\varepsilon\in(0,1]$ and $g\in C^\infty(\overline{Q_0},\mr^N)$. Then there exists a unique
    $$
    v\in C(\overline{I_0};L^2(B_0,\mr^N))\cap L^q(I_0;W^{1,q}(B_0,\mr^N))
    $$
    with $v-g\in L^q(I_0;W^{1,q}_0(B_0,\mr^N))$ and $v(\cdot,t_1)=g(\cdot,t_1)$ which is a weak solution to $v_t-\operatorname{div}\mathcal{A}_\varepsilon(t,Dv)=0$ in $Q_0$. Moreover, $Dv\in L^\infty_{\loc}(Q_0,\mr^{N\times n})$ and
    \begin{equation}\label{eq : estimate epsilon}
        \sup_{Q_{\rho/2}(z)}|Dv|\leq c\left(\miint{Q_\rho(z)}\left[\mathcal{G}_\varepsilon(t,|Dv|)|Dv|^2+1\right]dz\right)^{\frac12}
    \end{equation}
    for every $Q_\rho(z)\Subset Q_0$, where $c=c(n,q)$ is the constant in Proposition \ref{prop : a priori}.
\end{lemma}
\begin{proof}
The existence and uniqueness of $v$ follow from the theory of monotone operators \cite{Lions1969}. Indeed, the function $w\coloneq v-g$ must solve $w_t-\operatorname{div}\mathcal{A}_\varepsilon(t,Dg+Dw)=-g_t$ with $w(\cdot,t_1)=0$. The operator $w\mapsto-\operatorname{div}\mathcal{A}_\varepsilon(t,Dg+Dw)$ is monotone, hemicontinuous and bounded from $L^q(I_0;W^{1,q}_0(B_0,\mr^N))$ to its dual, and it is coercive because $a+\varepsilon\geq\varepsilon$. Uniqueness follows by testing the difference of two solutions with itself.

For $M\geq1$, let $v_M$ be the weak solution to $\partial_tv_M-\operatorname{div}\mathcal{A}_{\varepsilon,M}(t,Dv_M)=0$ in $Q_0$ with $v_M-g\in L^2(I_0;W^{1,2}_0(B_0,\mr^N))$ and $v_M(\cdot,t_1)=g(\cdot,t_1)$. It exists by the same theory in the Hilbert space setting, since $\mathcal{A}_{\varepsilon,M}(t,\cdot)$ is Lipschitz continuous and strongly monotone, as shown at the beginning of the proof of Proposition \ref{prop : a priori}. Let $L_g\coloneq\|Dg\|_{L^\infty(Q_0)}$. Then
$$
\begin{aligned}
|\mathcal{A}_{\varepsilon,M}(t,\xi)\cdot Dg|&\leq\frac12\mathcal{G}_{\varepsilon,M}(t,|\xi|)|\xi|^2\qquad\text{if } |\xi|\geq2L_g,\\
|\mathcal{A}_{\varepsilon,M}(t,\xi)\cdot Dg|&\leq2L_g^2\,\mathcal{G}_\varepsilon(t,2L_g)\qquad\text{if } |\xi|<2L_g.
\end{aligned}
$$
Testing the system for $v_M$ with $v_M-g$ and applying Gronwall's lemma to the term containing $g_t$, we obtain
\begin{equation}\label{eq : uniform energy}
\sup_{t\in I_0}\int_{B_0}|v_M-g|^2(x,t)\,dx+\iints{Q_0}\mathcal{G}_{\varepsilon,M}(t,|Dv_M|)|Dv_M|^2\,dz\leq C,
\end{equation}
with $C$ depending on $\varepsilon$, $g$, $Q_0$, $p$, $q$ and $\|a\|_{L^\infty(0,T)}$, but not on $M$. Proposition \ref{prop : a priori} applies to $v_M$, and together with \eqref{eq : uniform energy} it shows that for every $Q'\Subset Q_0$ there is a constant $K$, independent of $M$, such that $|Dv_M|\leq K$ in $Q'$. If $M\geq K$, then $\mathcal{G}_{\varepsilon,M}(t,|Dv_M|)=\mathcal{G}_\varepsilon(t,|Dv_M|)$ in $Q'$, so that $v_M$ solves $\partial_tv_M-\operatorname{div}\mathcal{A}_\varepsilon(t,Dv_M)=0$ in $Q'$. By \eqref{eq : monotonicity} and \eqref{eq : growth} with $\mu=\varepsilon$, the field $\mathcal{A}_\varepsilon(t,\cdot)$ is strongly monotone with constant $c\varepsilon^{p-2}$ and Lipschitz continuous on $\{|\xi|\leq K\}$ with a constant depending on $\varepsilon$, $K$, $p$, $q$ and $\|a\|_{L^\infty(0,T)}$. The difference quotient argument of Section \ref{section 3} therefore gives bounds for $D^2v_M$ in $L^2_{\loc}(Q')$, and consequently for $\partial_tv_M=\operatorname{div}\mathcal{A}_\varepsilon(t,Dv_M)$ in $L^2_{\loc}(Q')$, which are independent of $M$. Hence, for every cylinder $B\times J\Subset Q'$, the functions $Dv_M$ are bounded in $L^2(J;W^{1,2}(B,\mr^{N\times n}))$ and their time derivatives $\partial_tDv_M=D\partial_tv_M$ are bounded in $L^2(J;W^{-1,2}(B,\mr^{N\times n}))$.

By the compactness theorem of Simon \cite{Simon1987}, the functions $Dv_M$ are precompact in $L^2(B\times J)$, and a diagonal argument over an exhaustion of $Q_0$ gives a sequence $M_j\to\infty$ such that $Dv_{M_j}$ converges almost everywhere in $Q_0$, and $v_{M_j}-g$ converges weakly to some $\widetilde{v}-g$ in $L^2(I_0;W^{1,2}_0(B_0,\mr^N))$. Moreover, since $\mathcal{G}_{\varepsilon,M}(t,s)\leq\mathcal{G}_\varepsilon(t,s)\leq C(1+s)^{q-2}$, we have
$$
\left(\mathcal{G}_{\varepsilon,M}(t,s)s\right)^{q'}\leq C\left(\mathcal{G}_{\varepsilon,M}(t,s)s^2+1\right)
$$
with $C$ independent of $M$, and by \eqref{eq : uniform energy} the fluxes $\mathcal{A}_{\varepsilon,M_j}(t,Dv_{M_j})$ are bounded in $L^{q'}(Q_0)$. Since they converge almost everywhere to $\mathcal{A}_\varepsilon(t,D\widetilde{v})$, they converge weakly in $L^{q'}(Q_0)$. By Fatou's lemma and \eqref{eq : uniform energy}, $\iints{Q_0}\mathcal{G}_\varepsilon(t,|D\widetilde{v}|)|D\widetilde{v}|^2\,dz\leq C$, so that $D\widetilde{v}\in L^q(Q_0)$ and $\widetilde{v}-g\in L^q(I_0;W^{1,q}_0(B_0,\mr^N))$. Passing to the limit in the weak formulation for $v_{M_j}-g$, including test functions that do not vanish at $t=t_1$, we see that $\widetilde{v}$ is a weak solution to $\widetilde{v}_t-\operatorname{div}\mathcal{A}_\varepsilon(t,D\widetilde{v})=0$ in $Q_0$. Since the equation gives $\widetilde{v}_t\in L^{q'}(I_0;W^{-1,q'}(B_0,\mr^N))$, the function $\widetilde{v}$ belongs to $C(\overline{I_0};L^2(B_0,\mr^N))$ and satisfies $\widetilde{v}(\cdot,t_1)=g(\cdot,t_1)$. By uniqueness, $\widetilde{v}=v$. Finally, if $Q_\rho(z)\Subset Q_0$, then the functions $Dv_{M_j}$ are uniformly bounded on a neighborhood of $\overline{Q_\rho(z)}$ and converge almost everywhere. Thus we may pass to the limit in \eqref{eq : a priori estimate} for $v_{M_j}$, and we obtain \eqref{eq : estimate epsilon}.
\end{proof}

In the second approximation step we let $\varepsilon\to0$.
\begin{lemma}\label{lem : epsilon limit}
    There exist functions $g_\varepsilon\in C^\infty(\overline{Q_0},\mr^N)$, $\varepsilon\in(0,1]$, such that the solutions $u_\varepsilon$ given by Lemma \ref{lem : truncation limit} with $g=g_\varepsilon$ satisfy
    \begin{equation}\label{eq : strong convergence}
        \sup_{t\in I_0}\int_{B_0}|u_\varepsilon-u|^2(x,t)\,dx+\iints{Q_0}\left[|Du_\varepsilon-Du|^p+(a(t)+\varepsilon)|Du_\varepsilon-Du|^q\right]dz\longrightarrow0
    \end{equation}
    and $\mathcal{G}_\varepsilon(t,|Du_\varepsilon|)|Du_\varepsilon|^2\to H(t,|Du|)$ in $L^1(Q_0)$ as $\varepsilon\to0$.
\end{lemma}
\begin{proof}
Let $\phi_\delta$ be a standard mollifier in $\mr^{n+1}$ supported in $Q_\delta(0)$, and let $g^\delta\coloneq u*\phi_\delta$ on $\overline{Q_0}$ for $\delta>0$ so small that $Q_\delta(z)\Subset\Omega_T$ for every $z\in\overline{Q_0}$. Then $g^\delta\to u$ in $C(\overline{I_0};L^2(B_0,\mr^N))$ and in $L^q(I_0;W^{1,q}(B_0,\mr^N))$ as $\delta\to0$. Testing the weak formulation of \eqref{eq : main equation} with $\phi_\delta(z-\cdot)$ shows that
$$
\partial_tg^\delta=\operatorname{div}\left(\mathcal{A}(\cdot,Du)*\phi_\delta\right)\quad\text{in } Q_0,
$$
and we have
$$
f_\delta\coloneq\mathcal{A}(\cdot,Du)-\mathcal{A}(\cdot,Du)*\phi_\delta\to0\quad\text{in } L^{q'}(Q_0).
$$
We choose $\delta=\delta(\varepsilon)\in(0,\varepsilon]$ so small that $\varepsilon^{-\frac1{q-1}}\|f_{\delta(\varepsilon)}\|_{L^{q'}(Q_0)}^{q'}\leq\varepsilon$, and we set $g_\varepsilon\coloneq g^{\delta(\varepsilon)}$. In what follows we write $U_\varepsilon\coloneq Du_\varepsilon-Du$, $d_\varepsilon\coloneq Du-Dg_\varepsilon$, $f\coloneq f_{\delta(\varepsilon)}$ and
$$
\mathcal{E}\coloneq\iints{Q_0}\left[|U_\varepsilon|^p+(a(t)+\varepsilon)|U_\varepsilon|^q\right]dz.
$$

Since $u_\varepsilon-g_\varepsilon\in L^q(I_0;W^{1,q}_0(B_0,\mr^N))$ vanishes at $t=t_1$ and satisfies
$$
\partial_t(u_\varepsilon-g_\varepsilon)=\operatorname{div}\left(\mathcal{A}_\varepsilon(t,Du_\varepsilon)-\mathcal{A}(\cdot,Du)*\phi_{\delta(\varepsilon)}\right),
$$
testing with $u_\varepsilon-g_\varepsilon$ gives
$$
\frac12\int_{B_0}|u_\varepsilon-g_\varepsilon|^2(x,\tau)\,dx+\iints{Q_0^\tau}\left[\mathcal{A}_\varepsilon(t,Du_\varepsilon)-\mathcal{A}(t,Du)+f\right]\cdot(U_\varepsilon+d_\varepsilon)\,dz=0
$$
for every $\tau\in I_0$, where $Q_0^\tau\coloneq B_0\times(t_1,\tau)$. We split
$$
\mathcal{A}_\varepsilon(t,Du_\varepsilon)-\mathcal{A}(t,Du)=\left[\mathcal{A}_\varepsilon(t,Du_\varepsilon)-\mathcal{A}_\varepsilon(t,Du)\right]+r_\varepsilon,
$$
where $r_\varepsilon\coloneq\mathcal{A}_\varepsilon(t,Du)-\mathcal{A}(t,Du)$.
By \eqref{eq : monotonicity} with $\mu=\varepsilon$ and $m\in\{p,q\}$, and since $\varepsilon^2+|\xi|^2+|\eta|^2\geq\frac12|\xi-\eta|^2$, the first bracket satisfies
$$
\left[\mathcal{A}_\varepsilon(t,Du_\varepsilon)-\mathcal{A}_\varepsilon(t,Du)\right]\cdot U_\varepsilon\geq c\left(|U_\varepsilon|^p+(a(t)+\varepsilon)|U_\varepsilon|^q\right).
$$
The remaining terms are small. By \eqref{eq : growth} and $|Du_\varepsilon|\leq|Du|+|U_\varepsilon|$,
$$
\begin{aligned}
&\left|\mathcal{A}_\varepsilon(t,Du_\varepsilon)-\mathcal{A}_\varepsilon(t,Du)\right||d_\varepsilon|\\
&\qquad\leq C\left[(1+|Du|)^{p-2}|U_\varepsilon|+|U_\varepsilon|^{p-1}\right]|d_\varepsilon|\\
&\qquad\quad+C(a(t)+\varepsilon)\left[(1+|Du|)^{q-2}|U_\varepsilon|+|U_\varepsilon|^{q-1}\right]|d_\varepsilon|,
\end{aligned}
$$
and H\"older's and Young's inequalities give, for every $\vartheta\in(0,1)$ and with $p'\coloneq\frac{p}{p-1}$,
$$
\begin{aligned}
&\iints{Q_0}\left|\mathcal{A}_\varepsilon(t,Du_\varepsilon)-\mathcal{A}_\varepsilon(t,Du)\right||d_\varepsilon|\,dz\\
&\qquad\leq\vartheta\mathcal{E}+C_\vartheta\left(\|d_\varepsilon\|_{L^p(Q_0)}^{p'}+\|d_\varepsilon\|_{L^p(Q_0)}^p\right)\\
&\qquad\quad+C_\vartheta\left(\|d_\varepsilon\|_{L^q(Q_0)}^{q'}+\|d_\varepsilon\|_{L^q(Q_0)}^q\right),
\end{aligned}
$$
with $C_\vartheta$ depending also on $\|Du\|_{L^q(Q_0)}$ and $\|a\|_{L^\infty(0,T)}$. Next, we have
$$
|r_\varepsilon|\leq\varrho_p+a(t)\varrho_q+\varepsilon(\varepsilon^2+|Du|^2)^{\frac{q-1}{2}},\quad
\varrho_m\coloneq\left|(\varepsilon^2+|Du|^2)^{\frac{m-2}{2}}Du-|Du|^{m-2}Du\right|,
$$
where $\varrho_m$ tends to zero pointwise and is bounded by $C(1+|Du|)^{m-1}$. Young's inequality and the dominated convergence theorem yield
$$
\iints{Q_0}|r_\varepsilon||U_\varepsilon|\,dz\leq\vartheta\mathcal{E}+C_\vartheta\iints{Q_0}\left[\varrho_p^{p'}+a(t)\varrho_q^{q'}+\varepsilon(1+|Du|)^q\right]dz=\vartheta\mathcal{E}+o(1),
$$
and
$$
\iints{Q_0}|r_\varepsilon||d_\varepsilon|\,dz\leq\|r_\varepsilon\|_{L^{q'}(Q_0)}\|d_\varepsilon\|_{L^q(Q_0)}=o(1).
$$
Finally, since $\|U_\varepsilon\|_{L^q(Q_0)}\leq(\varepsilon^{-1}\mathcal{E})^{\frac1q}$,
$$
\begin{aligned}
\iints{Q_0}|f|(|U_\varepsilon|+|d_\varepsilon|)\,dz&\leq\vartheta\mathcal{E}+C_\vartheta\varepsilon^{-\frac1{q-1}}\|f\|_{L^{q'}(Q_0)}^{q'}+\|f\|_{L^{q'}(Q_0)}\|d_\varepsilon\|_{L^q(Q_0)}\\
&=\vartheta\mathcal{E}+o(1)
\end{aligned}
$$
by the choice of $\delta(\varepsilon)$. Choosing $\vartheta$ small and taking first $\tau=t_2$, we absorb the terms $\vartheta\mathcal{E}$ into the left-hand side and obtain $\mathcal{E}\to0$, since $\|d_\varepsilon\|_{L^q(Q_0)}\to0$. Returning to arbitrary $\tau$, we also obtain $\sup_\tau\|u_\varepsilon-g_\varepsilon\|_{L^2(B_0)}(\tau)\to0$, which together with $g_\varepsilon\to u$ in $C(\overline{I_0};L^2(B_0,\mr^N))$ proves \eqref{eq : strong convergence}.

For the convergence of the energy densities, we note that $\mathcal{E}\to0$ gives $Du_\varepsilon\to Du$ in $L^p(Q_0)$, $a^{\frac1q}Du_\varepsilon\to a^{\frac1q}Du$ in $L^q(Q_0)$ and $\varepsilon\iints{Q_0}|Du_\varepsilon|^q\,dz\to0$. Writing
$$
\mathcal{G}_\varepsilon(t,|Du_\varepsilon|)|Du_\varepsilon|^2=(\varepsilon^2+|Du_\varepsilon|^2)^{\frac{p-2}{2}}|Du_\varepsilon|^2+(a(t)+\varepsilon)(\varepsilon^2+|Du_\varepsilon|^2)^{\frac{q-2}{2}}|Du_\varepsilon|^2,
$$
we see that the first term converges to $|Du|^p$ and the $a(t)$-part of the second term converges to $a(t)|Du|^q$ in $L^1(Q_0)$ by Vitali's convergence theorem along every subsequence converging almost everywhere, while the $\varepsilon$-part tends to zero in $L^1(Q_0)$.
\end{proof}

We are now in a position to prove the Lipschitz bound.
\begin{proof}[Proof of Theorem \ref{thm : Lipschitz}]
Let $Q_r(z_0)\Subset\Omega_T$, and choose a cylinder $Q_0\Subset\Omega_T$ with $\overline{Q_r(z_0)}\subset Q_0$. Let $u_\varepsilon$ be the functions from Lemma \ref{lem : epsilon limit}. By Lemma \ref{lem : truncation limit},
$$
\sup_{Q_{r/2}(z_0)}|Du_\varepsilon|\leq c\left(\miint{Q_r(z_0)}\left[\mathcal{G}_\varepsilon(t,|Du_\varepsilon|)|Du_\varepsilon|^2+1\right]dz\right)^{\frac12}
$$
with $c=c(n,q)$. By Lemma \ref{lem : epsilon limit}, the right-hand side converges to that of \eqref{eq : Lipschitz estimate}, and a subsequence of $Du_\varepsilon$ converges to $Du$ almost everywhere in $Q_0$. This proves \eqref{eq : Lipschitz estimate}.
\end{proof}

Combining Theorems \ref{thm : second order} and \ref{thm : Lipschitz}, we find that the time derivative of a weak solution is locally square integrable, although the coefficient may be discontinuous in time.
\begin{corollary}\label{cor : time derivative}
    Let $a\in L^\infty(0,T)$ be nonnegative and let $u$ be a weak solution to \eqref{eq : main equation}. Then the flux $\mathcal{A}(\cdot,Du)$ has weak spatial derivatives in $L^2_{\loc}(\Omega_T)$, and $u_t=\operatorname{div}\mathcal{A}(\cdot,Du)\in L^2_{\loc}(\Omega_T,\mr^N)$. More precisely, if $Q_{2r}(z_0)\Subset\Omega_T$ and $|Du|\leq M$ in $Q_r(z_0)$ for some $M>0$, then
    \begin{equation}\label{eq : time derivative}
        \iints{Q_r(z_0)}|u_t|^2\,dz\leq\frac{c\,K_M}{r^2}\iints{Q_{2r}(z_0)}\left[|Du|^2+H(t,|Du|)\right]dz,
    \end{equation}
    where $K_M\coloneq M^{p-2}+\|a\|_{L^\infty(0,T)}M^{q-2}$ and $c=c(n,p,q)$.
\end{corollary}
\begin{proof}
For $m\geq2$, the map $W_m(\xi)\coloneq|\xi|^{\frac{m-2}{m}}\xi$ belongs to $C^1(\mr^{N\times n},\mr^{N\times n})$ with $|DW_m(\xi)|\leq C(m)|\xi|^{\frac{m-2}{m}}$, and $W_m(V_m(\xi))=|\xi|^{m-2}\xi$. Let $Q_{2r}(z_0)\Subset\Omega_T$ with $|Du|\leq M$ in $Q_r(z_0)$. Since $|V_m(Du)|\leq M^{\frac m2}$ in $Q_r(z_0)$, the derivative of $W_m$ is bounded by $C(m)M^{\frac{m-2}{2}}$ on the range of $V_m(Du)$. By Theorem \ref{thm : second order} and the chain rule, the function $|Du|^{p-2}Du=W_p(V_p(Du))$ has weak spatial derivatives in $Q_r(z_0)$ with
$$
\left|D\left(|Du|^{p-2}Du\right)\right|\leq CM^{\frac{p-2}{2}}|DV_p(Du)|.
$$
Similarly, for almost every $t$ with $a(t)>0$, we have $V_q(Du(\cdot,t))\in W^{1,2}(B_r(x_0))$ and
$$
\left|a(t)D\left(|Du|^{q-2}Du\right)\right|\leq C\|a\|_{L^\infty(0,T)}^{\frac12}M^{\frac{q-2}{2}}\,a(t)^{\frac12}|DV_q(Du)|,
$$
while $a(t)|Du|^{q-2}Du$ vanishes on the time slices where $a(t)=0$. Squaring, integrating over $Q_r(z_0)$ and using \eqref{eq : second order estimate}, we obtain
$$
\begin{aligned}
\iints{Q_r(z_0)}|D\mathcal{A}(t,Du)|^2\,dz&\leq CK_M\iints{Q_r(z_0)}\left[|DV_p(Du)|^2+a(t)|DV_q(Du)|^2\right]dz\\
&\leq\frac{CK_M}{r^2}\iints{Q_{2r}(z_0)}\left[|Du|^2+H(t,|Du|)\right]dz.
\end{aligned}
$$
The weak formulation of \eqref{eq : main equation} shows that $u_t=\operatorname{div}\mathcal{A}(\cdot,Du)$ in $Q_r(z_0)$, and \eqref{eq : time derivative} follows from $|\operatorname{div}F|^2\leq n|DF|^2$. Finally, by Theorem \ref{thm : Lipschitz}, every compact subset of $\Omega_T$ is covered by finitely many cylinders $Q_r(z_0)$ with $Q_{2r}(z_0)\Subset\Omega_T$ on which $Du$ is bounded.
\end{proof}

\section{\bf H\"older continuity of the gradient}\label{section 5}
In this section we prove Theorem \ref{thm : Holder}. By Theorem \ref{thm : Lipschitz}, the gradient is locally bounded, and it remains to pass from a bounded gradient to a H\"older continuous gradient. To this end, we compare the solution with solutions of frozen systems normalized at the level of the gradient bound, and we combine the comparison estimate with the oscillation decay of the frozen systems through a Campanato type interpolation.

For $0\leq\theta\leq1$, $s\geq0$ and $\xi,\eta\in\mr^{N\times n}$, we set
\begin{equation}\label{def : frozen}
\begin{aligned}
H_\theta(s)&\coloneq\theta s^p+(1-\theta)s^q,\qquad A_\theta(\xi)\coloneq\theta|\xi|^{p-2}\xi+(1-\theta)|\xi|^{q-2}\xi,\\
D_\theta(\xi,\eta)&\coloneq\theta|V_p(\xi)-V_p(\eta)|^2+(1-\theta)|V_q(\xi)-V_q(\eta)|^2.
\end{aligned}
\end{equation}
By \eqref{eq : V consequences}, applied separately with $m=p$ and $m=q$,
\begin{equation}\label{eq : frozen monotonicity}
\left(A_\theta(\xi)-A_\theta(\eta)\right)\cdot(\xi-\eta)\geq c\,D_\theta(\xi,\eta),
\end{equation}
and by \eqref{eq : V equivalence},
\begin{equation}\label{eq : weighted coercivity}
D_\theta(\xi,\eta)\geq c\left\{\theta(|\xi|+|\eta|)^{p-2}+(1-\theta)(|\xi|+|\eta|)^{q-2}\right\}|\xi-\eta|^2,
\end{equation}
where $c=c(p,q)>0$ does not depend on $\theta$. Since $|\xi|+|\eta|\geq|\xi-\eta|$ and $p\geq2$, the latter inequality gives in particular
\begin{equation}\label{eq : p coercivity}
D_\theta(\xi,\eta)\geq c\,\theta|\xi-\eta|^p.
\end{equation}

The autonomous ingredient is the regularity theory of Ok, Scilla and Stroffolini \cite{OkScillaStroffolini2024} for parabolic systems with Uhlenbeck structure and Orlicz growth. For the fields $A_\theta$, the constants in this theory are uniform in the phase parameter $\theta$. We write $Q_1\coloneq B_1\times(-1,1)$.
\begin{lemma}\label{lem : autonomous}
    Let $0\leq\theta\leq1$, and let $w\in C([-1,1];L^2(B_1,\mr^N))$ with $H_\theta(|Dw|)\in L^1(Q_1)$ be a weak solution to $w_s-\operatorname{div}A_\theta(Dw)=0$ in $Q_1$ satisfying
    $$
    \miint{Q_1}H_\theta(|Dw|)\,dz\leq K_0.
    $$
    Then there exist constants $C_0=C_0(n,N,p,q,K_0)$ and $\beta_0=\beta_0(n,N,p,q)\in(0,1)$, independent of $\theta$, such that
    \begin{equation}\label{eq : autonomous}
        \sup_{Q_{1/2}}|Dw|\leq C_0\qquad\text{and}\qquad\osc_{Q_\rho}Dw\leq C_0\rho^{\beta_0}\quad\text{for } 0<\rho\leq\frac14,
    \end{equation}
    where $Q_\rho=Q_\rho(0)$.
\end{lemma}
\begin{proof}
We verify the structural assumptions of \cite{OkScillaStroffolini2024}, which deals with systems, and we check that all constants are independent of $\theta$. For $s>0$, let
$$
\begin{aligned}
&\phi_\theta(s)\coloneq\frac{\theta}{p}s^p+\frac{1-\theta}{q}s^q,\qquad
\phi_\theta'(s)=\theta s^{p-1}+(1-\theta)s^{q-1},\\
&\phi_\theta''(s)=\theta(p-1)s^{p-2}+(1-\theta)(q-1)s^{q-2},
\end{aligned}
$$
so that $A_\theta(\xi)=\phi_\theta'(|\xi|)\xi/|\xi|$ for $\xi\neq0$ and
$$
p-1\leq\frac{s\phi_\theta''(s)}{\phi_\theta'(s)}\leq q-1,\qquad \frac1q\leq\phi_\theta(1)\leq\frac1p,\qquad \frac1qH_\theta(s)\leq\phi_\theta(s)\leq\frac1pH_\theta(s).
$$
The first quotient is a weighted average of $p-1$ and $q-1$, and its bounds do not deteriorate when one of the phases vanishes. Moreover, if $\eta\neq0$ and $|\xi-\eta|\leq|\eta|/2$, then
$$
|DA_\theta(\xi)-DA_\theta(\eta)|\leq C(p,q)\frac{|\xi-\eta|}{|\eta|}\phi_\theta''(|\eta|).
$$
Indeed, for $m\in\{p,q\}$ the map
$$
DA_m(\zeta)=|\zeta|^{m-2}\left(\mathrm{Id}+(m-2)\frac{\zeta\otimes\zeta}{|\zeta|^2}\right)
$$
is smooth away from the origin and homogeneous of degree $m-2$. Hence
$$
|D^2A_m(\zeta)|\leq C(m)|\zeta|^{m-3}\qquad\text{for } \zeta\neq0,
$$
and since $|\zeta|\geq|\eta|/2$ for every $\zeta$ on the segment from $\eta$ to $\xi$, we obtain
$$
|DA_m(\xi)-DA_m(\eta)|\leq C(m)|\eta|^{m-3}|\xi-\eta|.
$$
The claim follows by combining the two cases with the weights $\theta$ and $1-\theta$ and using $|\eta|^{m-2}\leq(m-1)|\eta|^{m-2}$. This is the relative Hessian condition of \cite[Assumption 2.3]{OkScillaStroffolini2024} with exponent $1$ and a constant independent of $\theta$. The index conditions of \cite{OkScillaStroffolini2024}, namely $1<\underline p<2<\overline q$ and $\underline p>\frac{2n}{n+2}$, are met after replacing $p$ and $q$ by $\underline p\coloneq2-\frac1{n+2}$ and $\overline q\coloneq q+1$, for which the preceding bounds for $s\phi_\theta''/\phi_\theta'$ remain valid since $p\geq2$.

The results of \cite{OkScillaStroffolini2024} are stated for functions normalized by $\phi(1)=1$, and their constants depend only on $n$, $N$, $\underline p$, $\overline q$ and the constant in the relative Hessian condition. We therefore pass to $\psi_\theta\coloneq\phi_\theta/c_\theta$ with $c_\theta\coloneq\phi_\theta(1)\in[\frac1q,\frac1p]$, which satisfies the same conditions with the same constants, and to the function $\widetilde{w}(y,\sigma)\coloneq w(y,\sigma/c_\theta)$, which solves
$$
\widetilde{w}_\sigma-\operatorname{div}\left(\psi_\theta'(|D\widetilde{w}|)\frac{D\widetilde{w}}{|D\widetilde{w}|}\right)=0\qquad\text{in } \widetilde{Q}\coloneq B_1\times(-c_\theta,c_\theta)
$$
and satisfies
$$
\miint{\widetilde{Q}}\psi_\theta(|D\widetilde{w}|)\,dz=\frac1{c_\theta}\miint{Q_1}\phi_\theta(|Dw|)\,dz\leq\frac qpK_0.
$$
The price of this normalization is that the time length of $\widetilde{Q}$ may be as small as $\frac2q$, and we work with cylinders whose radii are of order $q^{-\frac12}$. The cylinders in \cite{OkScillaStroffolini2024} are backward, and we write $Q^-_r(y,\sigma)\coloneq B_r(y)\times(\sigma-r^2,\sigma]$. Let $r_0\coloneq\frac1{8\sqrt q}$. For $|y|\leq\frac34$ and $|\sigma|\leq\frac9{16}c_\theta$, the cylinder $Q^-_{2r_0}(y,\sigma)$ is compactly contained in $\widetilde{Q}$, because $4r_0^2\leq\frac{c_\theta}{16}$, and \cite[Corollary 5.3]{OkScillaStroffolini2024} gives
$$
\|D\widetilde{w}\|_{L^\infty(Q^-_{r_0}(y,\sigma))}\leq c\left(\miint{Q^-_{2r_0}(y,\sigma)}\psi_\theta(|D\widetilde{w}|)\,dz+1\right)^{\frac{2}{(n+2)\underline p-2n}}
$$
with $c=c(n,N,p,q)$ and $(n+2)\underline p-2n=3$. Since $|\widetilde{Q}|\leq c(n,q)|Q^-_{2r_0}(y,\sigma)|$, the right-hand side is bounded by a constant $M_0=M_0(n,N,p,q,K_0)\geq1$. The set $B_{3/4}\times(-\frac9{16}c_\theta,\frac9{16}c_\theta)$ is the image of $Q_{3/4}$ under the change of variables, and it is covered by the cylinders $Q^-_{r_0}(y,\sigma)$. Hence $\sup_{Q_{3/4}}|Dw|\leq M_0$.

Let $0<\rho<\frac1{4\sqrt q}$ and $R_0\coloneq\frac1{2\sqrt q}$. The image of $Q_\rho$ is contained in $Q^-_{2\rho}(0,c_\theta\rho^2)$, since $2c_\theta\leq4$, and $Q^-_{R_0}(0,c_\theta\rho^2)$ is contained in the image of $Q_{3/4}$, since $R_0\leq\frac34$ and $R_0^2=\frac1{4q}\leq\frac{c_\theta}4$. Applying \cite[Theorem 6.1]{OkScillaStroffolini2024} in $Q^-_{R_0}(0,c_\theta\rho^2)$ with the gradient level $M\coloneq2M_0$, we obtain an exponent $\beta_0\in(0,1)$ depending only on $n,N,p,q$ and
$$
\begin{aligned}
\osc_{Q_\rho}Dw&\leq CM\left(\max\left\{\psi_\theta''(M)^{\frac12},\psi_\theta''(M)^{-\frac12}\right\}\frac{2\rho}{R_0}\right)^{\beta_0}\\
&\leq CM\max\left\{\psi_\theta''(M)^{\frac12},\psi_\theta''(M)^{-\frac12}\right\}\left(4\sqrt q\,\rho\right)^{\beta_0},
\end{aligned}
$$
where $C=C(n,N,p,q)$, and where we used that the maximum is at least one and $\beta_0<1$. Since $M\geq1$ and $c_\theta\in[\frac1q,\frac1p]$, we have $p\leq\psi_\theta''(M)\leq q(q-1)M^{q-2}$, so that the prefactor is bounded uniformly in $\theta$. For $\frac1{4\sqrt q}\leq\rho\leq\frac14$, we simply use $\osc_{Q_\rho}Dw\leq2M_0\leq2M_0\left(4\sqrt q\,\rho\right)^{\beta_0}$. This proves \eqref{eq : autonomous}.
\end{proof}

We now compare a weak solution whose gradient is bounded with solutions of frozen systems. Let $M\geq1$. For $z_*=(x_*,t_*)\in\Omega_T$ and $r>0$, we set
$$
g_*\coloneq M^{p-2}+a(t_*)M^{q-2},\qquad \theta_*\coloneq\frac{M^{p-2}}{g_*}\in(0,1],
$$
and
\begin{equation}\label{def : G cylinder}
G_r(z_*)\coloneq B_r(x_*)\times\left(t_*-\frac{r^2}{g_*},t_*+\frac{r^2}{g_*}\right).
\end{equation}
Since $g_*\geq1$, we have $G_r(z_*)\subset Q_r(z_*)$.
\begin{lemma}\label{lem : comparison}
    Let $a$ satisfy \eqref{cond : Holder in time} and let $u$ be a weak solution to \eqref{eq : main equation}. Suppose that $G_r(z_*)\Subset\Omega_T$, $|Du|\leq M$ in $G_r(z_*)$ and
    $$
    \omega\coloneq[a]_\beta M^{q-p}r^{2\beta}\leq1.
    $$
    Then there exists a weak solution $w$ to $w_t-\operatorname{div}\mathcal{A}(t_*,Dw)=0$ in $G_r(z_*)$ with $w=u$ on the parabolic boundary of $G_r(z_*)$ such that
    \begin{equation}\label{eq : comparison}
        \miint{G_r(z_*)}|Du-Dw|^p\,dz\leq CM^p\omega^2
    \end{equation}
    and
    \begin{equation}\label{eq : comparison oscillation}
        \osc_{G_\rho(z_*)}Dw\leq C_1M\left(\frac\rho r\right)^{\beta_0}\qquad\text{for } 0<\rho\leq\frac r4,
    \end{equation}
    where $C=C(p,q)$ and $C_1=C_1(n,N,p,q)$.
\end{lemma}
\begin{proof}
We normalize. For $(y,s)\in Q_1$, let $x=x_*+ry$, $t=t_*+r^2s/g_*$ and $\widetilde{u}(y,s)\coloneq u(x,t)/(Mr)$, so that $D\widetilde{u}(y,s)=Du(x,t)/M$ and $|D\widetilde{u}|\leq1$ in $Q_1$. A change of variables in the weak formulation shows that
$$
\widetilde{u}_s-\operatorname{div}\mathbb{A}(s,D\widetilde{u})=0\quad\text{in } Q_1,\qquad\mathbb{A}(s,\xi)\coloneq\theta_*|\xi|^{p-2}\xi+b(s)|\xi|^{q-2}\xi,
$$
where $b(s)\coloneq a(t_*+r^2s/g_*)M^{q-2}/g_*$ satisfies $b(0)=1-\theta_*$. Since $M^{q-2}/g_*=\theta_*M^{q-p}$ and $g_*\geq1$, \eqref{cond : Holder in time} gives
\begin{equation}\label{eq : theta factor}
|b(s)-(1-\theta_*)|\leq\theta_*M^{q-p}[a]_\beta\left(\frac{r^2|s|}{g_*}\right)^\beta\leq\theta_*\,\omega\qquad\text{for } |s|<1.
\end{equation}
Let $m=q$ if $\theta_*<1$ and $m=p$ if $\theta_*=1$. Since $|D\widetilde{u}|\leq1$, the function $\widetilde{u}$ belongs to $L^m(-1,1;W^{1,m}(B_1,\mr^N))$, and $\widetilde{u}_s$ belongs to $L^{m'}(-1,1;W^{-1,m'}(B_1,\mr^N))$. By the theory of monotone operators \cite{Lions1969}, there exists a unique $\widetilde{w}\in C([-1,1];L^2(B_1,\mr^N))$ with
$$
\widetilde{w}-\widetilde{u}\in L^m(-1,1;W^{1,m}_0(B_1,\mr^N))\quad\text{and}\quad \widetilde{w}(\cdot,-1)=\widetilde{u}(\cdot,-1)
$$
which solves $\widetilde{w}_s-\operatorname{div}A_{\theta_*}(D\widetilde{w})=0$ in $Q_1$. Testing the difference of the two systems with $\widetilde{u}-\widetilde{w}$ and using \eqref{eq : frozen monotonicity}, we obtain for every $\tau\in(-1,1)$
$$
\begin{aligned}
&\frac12\int_{B_1}|\widetilde{u}-\widetilde{w}|^2(y,\tau)\,dy+c\int_{-1}^\tau\int_{B_1}D_{\theta_*}(D\widetilde{u},D\widetilde{w})\,dy\,ds\\
&\qquad\leq\iints{Q_1}\left|A_{\theta_*}(D\widetilde{u})-\mathbb{A}(s,D\widetilde{u})\right||D\widetilde{u}-D\widetilde{w}|\,dy\,ds.
\end{aligned}
$$
By \eqref{eq : theta factor}, the integrand on the right-hand side is at most $\theta_*\omega|\xi|^{q-1}|\xi-\eta|$ with $\xi=D\widetilde{u}$ and $\eta=D\widetilde{w}$. If $\xi\neq0$, Young's inequality and $|\xi|+|\eta|\geq|\xi|$ give
$$
\theta_*\omega|\xi|^{q-1}|\xi-\eta|\leq\vartheta\theta_*(|\xi|+|\eta|)^{p-2}|\xi-\eta|^2+C_\vartheta\theta_*\omega^2|\xi|^{2q-p},
$$
and the same bound holds trivially if $\xi=0$. Letting $\tau\to1$, choosing $\vartheta=\vartheta(p,q)$ small, absorbing the first term by \eqref{eq : weighted coercivity}, and using $|\xi|\leq1$ and $2q-p\geq2$, we arrive at
\begin{equation}\label{eq : theta comparison}
\miint{Q_1}D_{\theta_*}(D\widetilde{u},D\widetilde{w})\,dy\,ds\leq C\theta_*\omega^2.
\end{equation}
Since $\theta_*>0$, the factor $\theta_*$ cancels against the factor $\theta_*$ in \eqref{eq : p coercivity}, and we obtain
$$
\miint{Q_1}|D\widetilde{u}-D\widetilde{w}|^p\,dy\,ds\leq C\omega^2.
$$
We now set $w(x,t)\coloneq Mr\widetilde{w}(y,s)$. Then $w$ solves $w_t-\operatorname{div}\mathcal{A}(t_*,Dw)=0$ in $G_r(z_*)$, since $1-\theta_*=a(t_*)M^{q-2}/g_*$, and $w=u$ on the parabolic boundary of $G_r(z_*)$. Changing variables back proves \eqref{eq : comparison}.

For \eqref{eq : comparison oscillation}, we note that $|V_m(\eta)|^2\leq2|V_m(\xi)|^2+2|V_m(\xi)-V_m(\eta)|^2$ for $m\in\{p,q\}$. Together with $|D\widetilde{u}|\leq1$ and \eqref{eq : theta comparison}, this gives
$$
\begin{aligned}
\miint{Q_1}H_{\theta_*}(|D\widetilde{w}|)\,dy\,ds&\leq2\miint{Q_1}H_{\theta_*}(|D\widetilde{u}|)\,dy\,ds+2\miint{Q_1}D_{\theta_*}(D\widetilde{u},D\widetilde{w})\,dy\,ds\\
&\leq2+2C.
\end{aligned}
$$
Lemma \ref{lem : autonomous} with $K_0=2+2C$ yields $\osc_{Q_\sigma}D\widetilde{w}\leq C_0\sigma^{\beta_0}$ for $0<\sigma\leq\frac14$. Since the image of $Q_\sigma$ under the change of variables is $G_{\sigma r}(z_*)$, this is \eqref{eq : comparison oscillation}.
\end{proof}

We are now ready to prove the main theorem.
\begin{proof}[Proof of Theorem \ref{thm : Holder}]
Let $Q_{4R}(z_0)\Subset\Omega_T$. Since $Q_{2R}(z)\subset Q_{4R}(z_0)$ for every $z\in Q_{2R}(z_0)$, Theorem \ref{thm : Lipschitz} gives
$$
M\coloneq\max\left\{1,\sup_{Q_{2R}(z_0)}|Du|\right\}\leq c\left(2^{n+2}\miint{Q_{4R}(z_0)}\left[H(t,|Du|)+1\right]dz\right)^{\frac12}.
$$
In the rest of the proof, $C$ denotes a constant depending only on $n$, $N$, $p$, $q$, $\beta$, $[a]_\beta$, $\|a\|_{L^\infty(0,T)}$, $R$ and $M$. We set
$$
R_b\coloneq\min\left\{R,\left([a]_\beta M^{q-p}\right)^{-\frac1{2\beta}}\right\},
$$
with the convention $R_b=R$ if $[a]_\beta=0$, and we fix $z_*\in Q_R(z_0)$. For $0<r\leq R_b$, we have $G_r(z_*)\subset Q_{2R}(z_0)$, and $\omega=[a]_\beta M^{q-p}r^{2\beta}\leq(r/R_b)^{2\beta}\leq1$. Let $w_r$ be the function given by Lemma \ref{lem : comparison}. By H\"older's inequality, $p\geq2$, \eqref{eq : comparison} and $\omega\leq(r/R_b)^{2\beta}$,
\begin{equation}\label{eq : L2 comparison}
\miint{G_r(z_*)}|Du-Dw_r|^2\,dz\leq\left(\miint{G_r(z_*)}|Du-Dw_r|^p\,dz\right)^{\frac2p}\leq CM^2\left(\frac r{R_b}\right)^{\frac{8\beta}{p}}.
\end{equation}
Let $0<\rho\leq r/4$. Using $(Dw_r)_{G_\rho(z_*)}$ as a comparison constant, together with \eqref{eq : L2 comparison}, \eqref{eq : comparison oscillation} and $|G_r(z_*)|/|G_\rho(z_*)|=(r/\rho)^{n+2}$, we obtain
\begin{equation}\label{eq : two scales}
\begin{aligned}
\miint{G_\rho(z_*)}|Du-(Du)_{G_\rho(z_*)}|^2\,dz
&\leq2\left(\frac r\rho\right)^{n+2}\miint{G_r(z_*)}|Du-Dw_r|^2\,dz\\
&\quad+2\left(\osc_{G_\rho(z_*)}Dw_r\right)^2\\
&\leq CM^2\left[\left(\frac r\rho\right)^{n+2}\left(\frac r{R_b}\right)^{\frac{8\beta}p}+\left(\frac\rho r\right)^{2\beta_0}\right].
\end{aligned}
\end{equation}
We put $\lambda\coloneq\rho/R_b$ and choose
$$
r=R_b\lambda^{\kappa_*},\qquad\kappa_*\coloneq\frac{n+2+2\beta_0}{n+2+2\beta_0+\frac{8\beta}p}\in(0,1).
$$
For $\lambda\leq4^{-\frac1{1-\kappa_*}}$, this choice satisfies $4\rho\leq r\leq R_b$, and both terms on the right-hand side of \eqref{eq : two scales} become multiples of $\lambda^{2\gamma}$, since
$$
2\beta_0(1-\kappa_*)=\frac{16\beta\beta_0}{p(n+2+2\beta_0)+8\beta}=2\gamma
\quad\text{and}\quad
(n+2)(\kappa_*-1)+\frac{8\beta}p\kappa_*=2\gamma.
$$
Hence
$$
\miint{G_\rho(z_*)}|Du-(Du)_{G_\rho(z_*)}|^2\,dz\leq CM^2\left(\frac\rho{R_b}\right)^{2\gamma}
$$
for $0<\rho\leq4^{-\frac1{1-\kappa_*}}R_b$, and, after enlarging $C$, also for $\rho\leq R_b$, since $|Du|\leq M$.

Finally, we return to standard parabolic cylinders. Let
$$
g_{\max}\coloneq M^{p-2}+\|a\|_{L^\infty(0,T)}M^{q-2}.
$$
For $h>0$ with $\sqrt{g_*}h\leq R_b$, we have $Q_h(z_*)\subset G_{\sqrt{g_*}h}(z_*)$ and
$$
\frac{|G_{\sqrt{g_*}h}(z_*)|}{|Q_h(z_*)|}=g_*^{\frac n2}\leq g_{\max}^{\frac n2}.
$$
Using the mean over the larger cylinder as a comparison constant, we find
$$
\begin{aligned}
\miint{Q_h(z_*)}|Du-(Du)_{Q_h(z_*)}|^2\,dz&\leq g_{\max}^{\frac n2}\miint{G_{\sqrt{g_*}h}(z_*)}|Du-(Du)_{G_{\sqrt{g_*}h}(z_*)}|^2\,dz\\
&\leq Ch^{2\gamma}
\end{aligned}
$$
for all $z_*\in Q_R(z_0)$ and $0<h\leq R_b/\sqrt{g_{\max}}$. The parabolic version of Campanato's characterization of H\"older continuous functions \cite{DaPrato1965} yields a representative of $Du$ in $C^{\gamma,\frac\gamma2}(Q_R(z_0))$ whose H\"older seminorm is bounded by $C$. Since $M$ is controlled by the energy on $Q_{4R}(z_0)$, this completes the proof.
\end{proof}

We conclude with a remark on continuous coefficients.
\begin{remark}\label{rem : VMO}
    If $a$ is merely continuous on $[0,T]$ with modulus of continuity $\omega_a$, the proofs of Lemma \ref{lem : comparison} and of Theorem \ref{thm : Holder} work with $\omega$ replaced by $M^{q-p}\omega_a(r^2)$, as long as this quantity is at most one. With $\rho=\tau r$ and $0<\tau\leq\frac14$, the estimate \eqref{eq : two scales} then becomes
    $$
    \miint{G_{\tau r}(z_*)}|Du-(Du)_{G_{\tau r}(z_*)}|^2\,dz\leq CM^2\left[\tau^{-(n+2)}\left(M^{q-p}\omega_a(r^2)\right)^{\frac4p}+\tau^{2\beta_0}\right],
    $$
    where $C$ does not depend on $z_*$. Choosing first $\tau$ small and then $r$ small, we find that the mean oscillation of $Du$ on $G_\rho(z_*)$ tends to zero as $\rho\to0$, uniformly in $z_*$. The passage to standard parabolic cylinders follows as in the last step of the proof of Theorem \ref{thm : Holder}, so that $Du$ has vanishing mean oscillation locally in $\Omega_T$. Continuity of $Du$ under a Dini type condition on $a$ requires a different iteration, which we do not pursue here.
\end{remark}

\vspace{0.5cm}
\noindent\textbf{Acknowledgments.} This work was supported by the National Research Foundation of Korea (NRF) grant funded by the Korea government [Grant Nos. RS-2025-00555316 and RS-2025-25415411].

\end{document}